\documentclass[11pt]{article}
\usepackage{amssymb,amsmath}
\usepackage{amscd}
\usepackage{mathrsfs}
\allowdisplaybreaks[1]
\newcommand{\cL}{\mathscr{L}}
\newcommand{\LL}{\mathcal{L}}
\newcommand{\C}{\mathbb{C}}
\newcommand{\N}{\mathbb{N}}
\newcommand{\R}{\mathbb{R}}
\newcommand{\U}{\mathbb{U}}
\newcommand{\V}{\mathbb{V}}
\newcommand{\X}{\mathbb{X}}
\newcommand{\dd}{\,{\rm d}}
\newcommand{\D}{{\rm d}}

\renewcommand{\Re}{\mathop{\mathrm{Re}}}

\numberwithin{equation}{section}
\newtheorem{thm}{Theorem}[section]

\newtheorem{prop}[thm]{Proposition}
\newtheorem{lem}[thm]{Lemma}

\newtheorem{cor}[thm]{Corollary}

\newcommand{\QED}{\mbox{}\hfill$\Box$}

\begin{document}

\title{Three-dimensional stability of Burgers vortices}

\author{
\null\\
Thierry Gallay \\ 
Institut Fourier \\
Universit\'e de Grenoble I \\
BP 74\\
38402 Saint-Martin-d'H\`eres, France \\
{\tt Thierry.Gallay@ujf-grenoble.fr}
\and
\\
Yasunori Maekawa\\    
Faculty of Science\\
Kobe University\\
1-1 Rokkodai, Nada-ku\\          
Kobe 657-8501, Japan\\
{\tt yasunori@math.kobe-u.ac.jp}}

\date{February 9, 2010}

\maketitle

\begin{abstract}
Burgers vortices are explicit stationary solutions of the Navier-Stokes
equations which are often used to describe the vortex tubes observed 
in numerical simulations of three-dimensional turbulence. In this 
model, the velocity field is a two-dimensional perturbation of a
linear straining flow with axial symmetry. The only free parameter 
is the Reynolds number $\rm{Re} = \Gamma/\nu$, where $\Gamma$ 
is the total circulation of the vortex and $\nu$ is the kinematic
viscosity. The purpose of this paper is to show that Burgers vortex
is asymptotically stable with respect to general three-dimensional 
perturbations, for all values of the Reynolds number. This definitive
result subsumes earlier studies by various authors, which were either
restricted to small Reynolds numbers or to two-dimensional 
perturbations. Our proof relies on the crucial observation that
the linearized operator at Burgers vortex has a simple and 
very specific dependence upon the axial variable. This allows to 
reduce the full linearized equations to a vectorial two-dimensional
problem, which can be treated using an extension of the techniques
developped in earlier works. Although Burgers vortices are found
to be stable for all Reynolds numbers, the proof indicates that
perturbations may undergo an important transient amplification if 
$\rm{Re}$ is large, a phenomenon that was indeed observed in numerical 
simulations. 
\end{abstract}

\section{Introduction}\label{sec.intro}

The axisymmetric Burgers vortex is an explicit solution of the
three-dimensional Navier-Stokes equations which provides a simple and
widely used model for the vortex tubes or filaments that are observed
in turbulent flows \cite{Bu,To}. Despite obvious limitations, due to
oversimplified assumptions, this model describes in a correct way the
fundamental mechanisms which are responsible for the persistence of
coherent structures in three-dimensional turbulence, namely the
balance between vorticity amplification due to stretching and
vorticity dissipation due to viscosity. If one believes that vortex
tubes play a significant role in the dynamics of turbulent flows, it
is an important issue to determine their stability with respect to
perturbations in the largest possible class. So far, this problem has
been studied only for the axisymmetric Burgers vortex and for a
closely related family of asymmetric vortices \cite{RS,MKO}.

As was shown by Leibovich and Holmes \cite{LH}, one cannot hope to
prove energetic stability of the Burgers vortex even if the 
circulation Reynolds number is very small. To tackle the stability
problem, it is therefore necessary to have a closer look at the
spectrum of the linearized operator. This is a relatively easy task if
we restrict ourselves to {\em two-dimensional} perturbations.
Assuming that the vortex tube is aligned with the vertical axis, this
means that the perturbed velocity field lies in the horizontal plane 
and does not depend on the vertical variable. Under such conditions, the
Burgers vortex is known to be stable for any value of the Reynolds
number. This result was first established by Giga and Kambe \cite{GK}
for $\rm{Re} \ll 1$ and then by Gallay and Wayne \cite{GW1} in the
general case.  Moreover, a lot is known about the spectrum of the
linearized operator, which turns out to be purely discrete in a
neighborhood of the origin in the complex plane. Using perturbative
expansions, Robinson and Saffman \cite{RS} showed that all linear
modes are exponentially damped for small Reynolds numbers. This property
was then numerically verified by Prochazka and Pullin \cite{PP1}
for $\rm{Re} \le 10^4$, and finally rigorously established in
\cite{GW1}.

The situation is much more complicated if we allow for arbitrary 
{\em three-dimensional} perturbations. In that case, it was shown by
Rossi and Le Diz\`es \cite{RD} that the linearized operator does not
have any eigenfunction with nontrivial dependence in the vertical
variable. While this result precludes the existence of unstable
eigenvalues, it also implies that stability cannot be deduced from
such a simple analysis, and that continuous spectrum necessarily plays
an important role. Unfortunately, the vertical dependence of the
perturbed solutions is not easy to determine, as can be seen from the
note \cite{Cr} where a few attempts are made in that direction. The
only rigorous result so far is due to Gallay and Wayne \cite{GW2}, who
proved that the Burgers vortex is asymptotically stable with respect
to three-dimensional perturbations in a fairly large class provided
that the Reynolds number is sufficiently small. For larger Reynolds
numbers, up to $\rm{Re} = 5000$, an important numerical work by
Schmid and Rossi \cite{SR} indicates that all modes are exponentially 
damped by the linearized evolution, although significant short-time 
amplification can occur.

In this paper, we prove that the axisymmetric Burgers vortex is
asymptotically stable with respect to {\em three-dimensional}
perturbations for {\em arbitrary values} of the Reynolds number. As in
\cite{GW2}, we assume that the perturbations are nicely localized in
the horizontal variables, but we do not impose any decay with respect 
to the vertical variable.  Our approach is based on the fact that the
linearized operator has a very simple dependence upon the vertical
variable: the only term involving $x_3$ is the dilation operator
$x_3 \partial_{x_3}$, which originates from the background straining
field. This crucial property was already exploited in \cite{RD,Cr,SR},
but we shall show that it allows to reduce the three-dimensional
stability problem to a two-dimensional one, which can then be treated
using an extension of the techniques developped in
\cite{GW1}. Although the spectrum of the linearized operator remains
stable for all Reynolds numbers, the estimates we have on the
associated semigroup deteriorate as $\rm{Re}$ increases, in
full agreement with the amplification phenomena observed in \cite{SR}.

We now formulate our results in a more precise way. We start from the 
three-dimensional incompressible Navier-Stokes equations:
\begin{equation}
  \partial_t V + (V,\nabla )V \,=\, \nu\Delta V - \frac{1}{\rho} \nabla P~, 
  \qquad \nabla \cdot V \,=\, 0~, \label{eq.NS}
\end{equation}
where $V = V(x,t) \in \R^3$ denotes the velocity field, $P = P(x,t)
\in \R$ is the pressure field, and $x = (x_1,x_2,x_3)^\top \in \R^3$
is the space variable. The parameters in \eqref{eq.NS} are the
kinematic viscosity $\nu > 0$ and the density $\rho > 0$. To obtain
tubular vortices, we assume that the velocity $V$ can be decomposed as
follows:
\begin{equation}
  V(x,t) \,=\, V^s(x) + U(x,t)~, \label{eq.velocity}
\end{equation}
where $V^s$ is an axisymmetric straining flow given by the explicit
formula
\begin{equation}
  V^s(x) \,=\, \frac{\gamma}{2} \begin{pmatrix} -x_1 \\ -x_2 \\ 2x_3
  \end{pmatrix} \,\equiv\, \gamma Mx~, \qquad \hbox{where} \quad
  M \,=\, \begin{pmatrix} -\frac12 & 0 & 0 \\ 0 & -\frac12 & 0 \\
  0 & 0 & 1\end{pmatrix}~.
  \label{def.strain}
\end{equation}
Here $\gamma > 0$ is a parameter which measures the intensity of the
strain. Note that $\nabla\cdot V^s = 0$, and that $V^s$ is a
stationary solution of \eqref{eq.NS} with the associated pressure $P^s
= -\frac12 \rho |V^s|^2$. Our goal is to study the evolution of the
perturbed velocity field $U(x,t)$.

To simplify the notations, we shall assume henceforth that 
$\gamma = \nu = \rho = 1$. This can be achieved without loss of 
generality by replacing the variables $x$, $t$ and the functions 
$V$, $P$ with the dimensionless quantities
\[
  \tilde x \,=\, \Bigl(\frac{\gamma}{\nu}\Bigr)^{1/2}x~, \qquad 
  \tilde t \,=\, \gamma t~, \qquad
  \tilde V \,=\, \frac{V}{(\gamma \nu)^{1/2}}~, \qquad 
  \tilde P \,=\, \frac{P}{\rho \gamma \nu}~.
\]
For further convenience, instead of considering the evolution of $V$ 
or $U$, we prefer working with the vorticity field $\Omega = 
\nabla \times V = \nabla \times U$. Taking the curl of \eqref{eq.NS} 
and using \eqref{eq.velocity}, \eqref{def.strain}, we obtain for 
$\Omega$ the evolution equation
\begin{equation}
  \partial_t \Omega + (U,\nabla )\Omega - (\Omega,\nabla )U \,=\, 
  L\Omega~, \qquad \nabla \cdot \Omega \,=\, 0~, \label{eq.vortex}
\end{equation} 
where $L$ is the differential operator defined by
\begin{equation}
  L\Omega \,=\, \Delta \Omega - (Mx,\nabla )\Omega + M\Omega~.
  \label{def.L}
\end{equation}

Under mild assumptions that will be specified below, the 
velocity field $U$ can be recovered from the vorticity $\Omega$ 
via the three-dimensional Biot-Savart law
\begin{equation}
  U(x) \,=\, -\frac{1}{4\pi}\int_{\R^3}\frac{(x-y)\times 
  \Omega(y)}{|x-y|^3} \dd y ~=:~ (K_{3D}*\Omega)(x)~.
  \label{def.3DBS}
\end{equation}
In what follows we shall often encounter the particular situation 
where the velocity $U$ is two-dimensional and horizontal, 
namely $U(x) = (U_1(x_h),U_2(x_h),0)^\top$ where $x_h = 
(x_1,x_2)^\top \in \R^2$. In that case the vorticity satisfies
$\Omega(x) = (0,0,\Omega_3(x_h))^\top$, and the relation 
\eqref{def.3DBS} reduces to the two-dimensional Biot-Savart law
\begin{equation}
  U_h(x_h) \,=\, \frac{1}{2\pi}\int_{\R^2} \frac{(x_h-y_h)^\bot}
  {|x_h-y_h|^2}\,\Omega_3(y_h) \dd y_h ~=:~ (K_{2D}\star\Omega_3)(x_h)~,
  \label{def.2DBS}
\end{equation}
where $U_h =(U_1,U_2)^\top$ and $x_h^\bot = (-x_2,x_1)^\top$. 

We can now introduce the {\em Burgers vortices}, which are explicit 
stationary solutions of \eqref{eq.vortex} of the form $\Omega = \alpha G$, 
where $\alpha \in \R$ is a parameter. The vortex profile is given
by 
\begin{equation}
  \quad~~ G(x) \,=\, \begin{pmatrix} 0 \\ 0 \\ g(x_h)\end{pmatrix}~, 
  \qquad \hbox{where} \quad
  g(x_h) \,=\, \frac{1}{4\pi}\,e^{-|x_h|^2/4}~. 
  \label{def.g}
\end{equation}
The associated velocity field $U = \alpha U^G$ can be obtained 
from the Biot-Savart law \eqref{def.2DBS} and has the following
form
\begin{equation}
  U^G(x) \,=\, u^g(|x_h|^2)\begin{pmatrix} -x_2 \\ x_1 \\ 0\end{pmatrix}~, 
  \qquad \hbox{where} \quad u^g(r) \,=\, \frac{1}{2\pi r}
  \Bigl(1-e^{-r/4}\Bigr)~. \label{def.u^g}
\end{equation}
If $\Omega = \alpha G$, it is easy to verify that $\alpha = 
\int_{\R^2}\Omega_3(x_h)\dd x_h$. This means that the parameter 
$\alpha \in \R$ represents the total circulation of the Burgers vortex 
$\alpha G$. In the physical literature, the quantity $|\alpha|$ is 
often referred to as the (circulation) Reynolds number. 

The aim of this paper is to study the asymptotic stability of the 
Burgers vortices. We thus consider solutions of \eqref{eq.vortex} 
of the form $\Omega = \alpha G + \omega$, $U = \alpha U^G + u$, 
and obtain the following evolution equation for the perturbation: 
\begin{equation}
  \partial_t\omega +(u,\nabla)\omega - (\omega,\nabla)u \,=\, 
  (L-\alpha \Lambda )\omega~, \qquad \nabla \cdot \omega \,=\,0~,
  \label{eq.omega}
\end{equation}
where $\Lambda$ is the integro-differential operator defined by
\begin{equation}
  \Lambda \omega \,=\, (U^G,\nabla )\omega - (\omega,\nabla )U^G + 
  (u,\nabla )G - (G,\nabla)u~. \label{def.Lambda}
\end{equation}
Here and in the sequel, it is always understood that $u =
K_{3D}*\omega$.  

An important issue is now to fix an appropriate function space for the
admissible perturbations. Since the Burgers vortex itself is
essentially a two-dimensional flow, it is natural to choose a
functional setting which allows for perturbations in the same class,
but we also want to consider more general ones.  Following \cite{GW2},
we thus assume that the perturbations are nicely localized in the
horizontal variables, but merely bounded in the vertical direction. As
we shall see below, this choice is more or less imposed by the
particular form of the linear operator \eqref{def.L}.

To specify the horizontal decay of the admissible perturbations, 
we first introduce two-dimensional spaces. Given $m \in [0,\infty]$, 
we denote by $\rho_m : [0,\infty) \to [1,\infty)$ the weight function
defined by 
\begin{equation}\label{def.weight}
  \rho_m(r) \,=\, \begin{cases} 1 & \hbox{if} \quad m = 0\,, \\
  (1+\frac{r}{4m})^m  & \hbox{if} \quad 0 < m<\infty\,, \\
  e^{r/4} & \hbox{if} \quad m = \infty\,.\end{cases}
\end{equation}
We introduce the weighted $L^2$ space 
\begin{equation}
  L^2(m) \,=\, \Bigl\{f\in L^2(\R^2)~\Big|~\int_{\R^2} |f(x_h)|^2
  \rho_m(|x_h|^2)\dd x_h < \infty \Bigr\}~, \label{def.L^2(m)}
\end{equation}
which is a Hilbert space with a natural inner product. Using 
H\"older's inequality, it is easy to verify that $L^2(m) 
\hookrightarrow L^1(\R^2)$ if $m > 1$. In that case, we also 
define the closed subspace
\begin{equation}
  L^2_0(m) \,=\, \Bigl\{f \in L^2(m)~\Big|~\int_{\R^2} f(x_h)
  \dd x_h = 0\Bigr\}~. \label{def.L^2_0(m)}
\end{equation}

Next, we define the three-dimensional space $X(m)$ as the set of all
$\phi : \R^3 \to \R$ for which the map $x_h \mapsto \phi(x_h,x_3)$
belongs to $L^2(m)$ for any $x_3 \in \R$, and is a bounded and
continuous function of $x_3$. In other words, we set
\begin{equation}
  X(m) \,=\, BC(\R\,;L^2(m))~, \qquad  
  X_0(m) \,=\, BC(\R\,;L^2_0(m))~, \label{def.X(m)}
\end{equation}
where ``$BC(\R\,;Y)$'' denotes the space of all bounded and continuous
functions from $\R$ into $Y$. Both $X(m)$ and $X_0(m)$ are Banach
spaces equipped with the norm
\begin{equation}\label{def.Xnorm}
  \|\phi\|_{X(m)} \,=\, \sup_{x_3\in\R}\|\phi(\cdot,x_3)\|_{L^2(m)}~.
\end{equation}

Our goal is to study the stability of the Burgers vortex $\Omega =
\alpha G$ with respect to perturbations $\omega \in X(m)^3$. In fact, 
we can assume without loss of generality that $\omega$ belongs to 
the subspace
\begin{equation}
  \X(m) \,=\, X(m) \times X(m) \times X_0(m) \,\subset\,
  X(m)^3~, \label{def.mathX(m)}
\end{equation}
which is invariant under the evolution defined by \eqref{eq.omega}. 
This is a consequence of the following result, whose proof is 
postponed to Section~\ref{subsec.lem.reduce}:

\begin{lem}\label{lem.reduction} Fix $m \in (1,\infty]$. If 
$\tilde\omega \in X(m)^3$ satisfies $\nabla \cdot \tilde\omega = 0$ 
in the sense of distributions, then there exists $\tilde \alpha \in \R$ 
such that
\begin{equation}\label{def.tildealpha}
  \int_{\R^2}\tilde\omega_3(x_h,x_3)\dd x_h \,=\, \tilde\alpha~,  
  \quad \hbox{for all } x_3\in \R~.
\end{equation}
\end{lem}

In view of Lemma~\ref{lem.reduction}, if $\Omega = \alpha G
+ \tilde\omega$ for some $\tilde\omega \in X(m)^3$, we can write
$\Omega = (\alpha +\tilde\alpha)G + \omega$, where $\tilde
\alpha$ is given by \eqref{def.tildealpha} and $\omega = \tilde 
\omega - \tilde\alpha G$. Then $\omega \in \X(m)$ by construction, 
and we are led back to the stability analysis of the Burgers 
vortex $(\alpha +\tilde\alpha)G$ with respect to perturbations 
in $\X(m)$.

In what follows we always consider the solutions $\omega(x,t)$ of
\eqref{eq.omega} as $\X(m)$-valued functions of time, and we often
denote by $\omega(\cdot,t)$ or simply $\omega(t)$ the map $x \mapsto
\omega(x,t)$. A minor drawback of our functional setting is that we
cannot expect the solutions of \eqref{eq.omega} to be continuous in
time in the strong topology of $\X(m)$. This is because the operator
$L$ defined in \eqref{def.L} contains the dilation operator
$-x_3\partial_{x_3}$, see Section~\ref{subsec.semigroup.L} below. To
restore continuity, it is thus necessary to equip $\X(m)$ with a
weaker topology. Following \cite{GW2}, we denote by $X_{loc}(m)$ the
space $X(m)$ equipped with the topology defined by the family of
seminorms
\[
  \|\phi\|_{X_n(m)} \,=\, \sup_{|x_3|\le n}\|\phi(\cdot,x_3)\|_{L^2(m)}~,
  \qquad n\in \N~. 
\]
In analogy with \eqref{def.mathX(m)}, we set ${\X}_{loc}(m) = X_{loc}(m) 
\times X_{loc}(m) \times X_{0,loc}(m)$, where $X_{0,loc}(m)$ is of course
the space $X_0(m)$ equipped with the topology of $X_{loc}(m)$. 

\medskip

We are now able to formulate our main result:

\begin{thm}\label{thm.main1} Fix $m \in (2,\infty]$ and $\alpha\in \R$. 
Then there exist $\delta = \delta(\alpha,m) > 0$ and $C = 
C(\alpha,m) \ge 1$ such that, for any $\omega_0 \in \X(m)$ with 
$\nabla \cdot \omega_0 = 0$ and $\|\omega_0\|_{\X(m)}\le \delta$, 
Eq.~\eqref{eq.omega} has a unique solution $\omega \in 
L^\infty(\R_+\,;\X(m))\cap C([0,\infty)\,;\X_{loc}(m))$ with initial 
data $\omega_0$. Moreover,
\begin{equation}
  \|\omega(t)\|_{\X(m)} \,\le\, C\|\omega_0\|_{\X(m)}\,e^{-t/2}~,
  \qquad \hbox{for all } t \ge 0~.
  \label{est.thm.main1}
\end{equation}
\end{thm}

Theorem~\ref{thm.main1} shows that the Burgers vortex $\alpha G$
is {\em asymptotically stable} with respect to perturbations 
in $\X(m)$, for any value of the circulation $\alpha \in \R$. 
If one prefers to consider perturbations in the larger space 
$X(m)^3$, then our result means that the family $\{\alpha G\}_{\alpha 
\in \R}$ of all Burgers vortices is asymptotically stable 
{\em with shift}, because the perturbations may then modify the 
circulation of the underlying vortex. The key point in the 
proof is to show that the linearized operator $L - \alpha\Lambda$ 
has a {\em uniform spectral gap} for all $\alpha \in \R$. This
implies a uniform decay rate in time for the perturbations, 
as in \eqref{est.thm.main1}. However, it should be emphasized
that the constants $C$ and $\delta$ in Theorem~\ref{thm.main1}
do depend on $\alpha$, in such a way that $C(\alpha,m) \to \infty$
and $\delta(\alpha,m) \to 0$ as $|\alpha| \to \infty$. This is 
in full agreement with the amplification phenomena numerically 
observed in \cite{SR}. 

The proof of Theorem~\ref{thm.main1} gives a more detailed information
on the solutions of \eqref{eq.omega} than what is summarized in
\eqref{est.thm.main1}. First of all, we can prove stability for any 
$m > 1$, but the exponential factor $e^{-t/2}$ in \eqref{est.thm.main1}
should then be replaced by $e^{-\eta t}$, where $\eta < (m-1)/2$ if $m
\le 2$. Next, thanks to parabolic smoothing, we can obtain decay
estimates not only for $\omega(t)$ but also for its spatial
derivatives. Finally, due to the particular structure of the linear
operator $L - \alpha\Lambda$, it turns out that the horizontal part
$\omega_h = (\omega_1,\omega_2)^\top$ of the vorticity vector has a
faster decay than the vertical component $\omega_3$ as $t \to
\infty$. Thus, a more complete (but less readable) version of our
result is as follows:

\begin{thm}\label{thm.main2} Fix $m \in (1,\infty]$, $\alpha\in \R$, 
and take $\mu \in (1,\frac32)$, $\eta \in (0,\frac12]$ such that 
$2\mu < m+1$ and $2\eta < m-1$. Then there exist $\delta = 
\delta(\alpha,m) > 0$ and $C = C(\alpha,m,\mu,\eta) > 1$ such that, 
for all initial data $\omega_0 \in \X(m)$ with $\nabla \cdot \omega_0 = 0$ 
and $\|\omega_0\|_{\X(m)}\le \delta$, Eq.~\eqref{eq.omega} has a unique 
solution $\omega \in L^\infty(\R_+\,;\X(m))\cap C([0,\infty)\,;\X_{loc}(m))$.  
Moreover, for all $t > 0$, 
\begin{align}
  \|\partial_x^\beta\omega_h(t)\|_{X(m)^2} \,&\le\, 
  \frac{C\|\omega_0\|_{\X(m)}}{a(t)^{|\beta|/2}} \,e^{-\mu t}~,
  \label{est.thm.main2.1}\\
  \|\partial_x^\beta\omega_3(t)\|_{X(m)} \,&\le\,
  \frac{C\|\omega_0\|_{\X(m)}}{a(t)^{|\beta|/2}}\,e^{-\eta t}~,
  \label{est.thm.main2.2}
\end{align}
where $a(t) = 1-e^{-t}$ and $\beta \in \N^3$ is any multi-index of 
length $|\beta| = \beta_1 + \beta_2 + \beta_3 \le 1$. 
\end{thm}

The decay rates \eqref{est.thm.main2.1}, \eqref{est.thm.main2.2} 
are optimal when $\beta = 0$, but it turns out that vertical 
derivatives such as $\partial_{x_3}\omega_h(t)$ or $\partial_{x_3}
\omega_3(t)$ have a faster decay as $t \to \infty$, see 
Sections~\ref{sec.linear} and \ref{sec.nonlinear} for 
more details. In any case, we believe that the optimal rates 
are those provided by the linear stability analysis, as in 
Proposition~\ref{prop.linear} below. 

The rest of this paper is devoted to the proof of
Theorems~\ref{thm.main1} and \ref{thm.main2}. Before giving the 
details, we explain here the main ideas in an informal way.
As was already mentioned, the main difficulty is to obtain 
good estimates on the solutions of the linearized equation
\begin{equation}
  \partial_t\omega \,=\, (L - \alpha\Lambda)\omega~, \qquad 
  \nabla \cdot \omega \,=\,0~.
  \label{eq.linear}
\end{equation}
Once this is done, the nonlinear terms in \eqref{eq.omega} can be
controlled using rather standard arguments, which are recalled in
Section~\ref{sec.nonlinear}. To study \eqref{eq.linear}, we use the
fact that the operator $L - \alpha\Lambda$ depends on the vertical
variable in a simple and very specific way.  Indeed, it is easy to
verify that $[\partial_{x_3},L] = -\partial_{x_3}$ and
$[\partial_{x_3},\Lambda] = 0$. This key observation, which already
plays a crucial role in the previous works \cite{RD,Cr,SR}, implies
the following identity:
\begin{equation}
  \partial_{x_3}^k \,e^{t(L-\alpha\Lambda)}\omega_0 \,=\, 
  e^{-kt}\,e^{t(L-\alpha\Lambda)}\partial_{x_3}^k \omega_0~,
  \label{eq.semigroup}
\end{equation}
for all $k \in \N$ and all $t \ge 0$. If we take $k \in \N$ sufficiently
large, depending on $|\alpha|$, we can use \eqref{eq.semigroup} to
to show that $\partial_{x_3}^k \omega(t)$ decays exponentially as 
$t \to \infty$ if $\omega(t)$ is a solution of \eqref{eq.linear}.
Then, by an interpolation argument, we deduce that all expressions
involving at least one vertical derivative play a negligible role
in the long-time asymptotics, see Section~\ref{sec.linear} for more 
details. This ``smoothing effect'' in the vertical direction is 
due to the stretching properties of the linear flow \eqref{eq.velocity}. 

As a consequence of these remarks, we can restrict our attention to
those solutions of \eqref{eq.linear} which are independent of the
vertical variable $x_3$. We call this particular situation the {\em
vectorial 2D problem}, and we study it in Section~\ref{sec.vectorial}. 
Note that the perturbations we consider here are two-dimensional in
the sense that $\partial_{x_3} u = \partial_{x_3} \omega = 0$, but
that all three components of $u$ or $\omega$ are possibly
nonzero. This is in contrast with the purely two-dimensional case
considered in \cite{GW1,GW2}, where in addition $u_3 = \omega_1 =
\omega_2 = 0$. Nevertheless, it is possible to show that the solutions
of \eqref{eq.linear} with $\partial_{x_3} \omega = 0$ converge
exponentially to zero as $t \to \infty$, and that the decay rate is
uniform in $\alpha$.  Extending the techniques developped in
\cite{GW1,GW2}, this can be done using spectral estimates and a
detailed study of the eigenvalue equation $(L-\alpha\Lambda)\omega =
\lambda\omega$. It is then a rather straightforward task to complete
the proof of Theorem~\ref{thm.main1} using the arguments presented
above.

\medskip\noindent{\bf Remark.} The vortex tubes observed in numerical
simulations are usually not axisymmetric: in general, they rather
exhibit an elliptical core region. A simple model for such asymmetric 
vortices is obtained by replacing the straining flow $V^s$ in 
\eqref{def.strain} with the nonsymmetric strain $V^s_\lambda(x) = 
\gamma M_\lambda x$, where $\lambda \in (0,1)$ is an asymmetry 
parameter and 
\begin{equation}
  M_\lambda \,=\, \begin{pmatrix} -\frac{1+\lambda}{2} & 0 & 0 \\ 0 & 
  -\frac{1-\lambda}{2} & 0 \\ 0 & 0 & 1\end{pmatrix}~.
  \label{def.strain'}
\end{equation}
Asymmetric Burgers vortices are then stationary solutions to 
\eqref{eq.vortex}, where the operator $L$ in the right-hand
side is defined by \eqref{def.L} with $M$ replaced by $M_\lambda$. 
Unlike in the symmetric case $\lambda = 0$, no explicit formula is 
available and proving the existence of stationary solutions is 
already a nontrivial task, except perhaps in the perturbative 
regime where either the asymmetry parameter $\lambda$ or the 
circulation number $\alpha$ is very small. In view of these difficulties, 
asymmetric Burgers vortices were first studied using formal asymptotic 
expansions and numerical calculations, see e.g. \cite{RS,MKO,PP2}. 
The mathematical theory is more recent, and includes several 
existence results which cover now the whole range of parameters 
$\lambda \in (0,1)$ and $\alpha \in \R$ \cite{GW2,GW3,M1,M2}. 
In addition, the stability with respect to two-dimensional perturbations 
is known to hold at least for small values of the asymmetry parameter 
\cite{GW3,M1}. However, the only result so far on three-dimensional 
stability is restricted to the particular case where the circulation 
number $\alpha$ is sufficiently small, depending on $\lambda$ 
\cite{GW2}. 

Using Theorem~\ref{thm.main1} and a simple perturbation argument, 
it is easy to show that asymmetric Burgers vortices are stable
with respect to three-dimensional pertubations in the space 
$\X(m)$, provided that the asymmetry parameter $\lambda$ is 
small enough depending on the circulation number $\alpha$. 
This follows from the fact the the linearized operator at the
symmetric Burgers vortex has a uniform spectral gap for 
all $\alpha \in \R$, and that the asymmetric Burgers vortex is 
$O(\lambda)$ close to the corresponding symmetric vortex in 
the topology of $\X(m)$, uniformly for all $\alpha \in \R$ \cite{GW3}. 
Although this stability result is new and not covered by \cite{GW2}, 
it is certainly not optimal, and we prefer to postpone the study 
of the three-dimensional stability of asymmetric Burgers vortices
to a future investigation. 


\section{Preliminaries}\label{sec.pre}

In this preliminary section we collect a few basic estimates which
will be used throughout the proof of Theorems~\ref{thm.main1} and
\ref{thm.main2}. They concern the semigroup generated by the linear
operator \eqref{def.L}, and the Biot-Savart law \eqref{def.3DBS}
relating the velocity field to the vorticity. Most of the results
were already established in \cite[Appendix~A]{GW2}, and are reproduced
here for the reader's convenience.

As in \cite{GW2}, we introduce the following generalization of 
the function spaces \eqref{def.L^2(m)} and \eqref{def.X(m)}. Given 
$m \in [0,\infty]$ and $p \in [1,\infty)$, we define the weighted
$L^p$ space
\[
  L^p(m) \,=\, \Bigl\{f\in L^p(\R^2)~\Big| \|f\|_{L^p(m)}^p = 
  \int_{\R^2}|f(x_h)|^p \rho_m(|x_h|^2)^{p/2} \dd x_h < \infty
  \Bigr\}~, 
\]
and the corresponding three-dimensional space
\[
  X^p(m) \,=\, BC(\R\,;L^p(m))~, \qquad  
  \|\phi\|_{X^p(m)} \,=\, \sup_{x_3\in\R}\|\phi(\cdot,x_3)\|_{L^p(m)}~.
\]
If $m > 2-\frac{2}{p}$, we also denote by $L^p_0(m)$ the subspace of all
$f \in L^p(m)$ such that $\int_{\R}f\dd x_h = 0$. In analogy 
with \eqref{def.mathX(m)}, we set $\X^p(m) \,=\, X^p(m) \times X^p(m) 
\times X_0^p(m)$ , where $X_0^p(m) = BC(\R\,;L^p_0(m))$. 

\subsection{The semigroup generated by $L$}\label{subsec.semigroup.L}

If we decompose the vorticity $\omega$ into its horizontal part
$\omega_h = (\omega_1,\omega_2)^\top$ and its vertical component 
$\omega_3$, it is clear from \eqref{def.strain} and \eqref{def.L}
that the linear operator $L$ has the following expression:
\begin{equation}
  L\omega \,=\, \begin{pmatrix} L_h\omega_h \\ L_3\omega_3\end{pmatrix}
  \,=\, \begin{pmatrix} (\LL_h + \LL_3 -\frac{3}{2})\omega_h \\
  (\LL_h + \LL_3)\omega_3\end{pmatrix}~, \label{eq.Lexp}
\end{equation}
where $\LL_h$ is the two-dimensional Fokker-Planck operator
\begin{equation}\label{def.LLh}
  \LL_h \,=\, \Delta_h + \frac{x_h}{2}\cdot\nabla_h + 1 \,=\, 
  \sum_{j=1}^2\partial_{x_j}^2 + \sum_{j=1}^2\frac{x_j}{2}\partial_{x_j}
  + 1~,
\end{equation}
and $\LL_3 = \partial_{x_3}^2 - x_3\partial_{x_3}$ is a convection-diffusion
operator in the vertical variable. 

As is shown in \cite[appendix A]{GW0}, the operator $\LL_h$ is
the generator of a strongly continuous semigroup in $L^2(m)$ given 
by the explicit formula
\begin{equation}
  (e^{t\LL_h}f)(x_h) \,=\, \frac{e^t}{4\pi a(t)} \int_{\R^2}
  e^{-\frac{|x_h-y_h|^2}{4a(t)}} f(y_he^{t/2})\dd y_h~, \qquad 
  t > 0~, \label{eq.semigroup.calL_h}
\end{equation}
where $a(t)=1-e^{-t}$. Similarly, the operator $\LL_3$ generates 
a semigroup of contractions in $BC(\R)$ given by
\begin{equation}
  (e^{t\LL_3}f)(x_3) \,=\, \frac{1}{\sqrt{2\pi a(2t)}}\int_{\R}
  e^{-\frac{|x_3 e^{-t}-y_3|^2}{2a(2t)}}f(y_3)\dd y_3~,  \qquad
  t > 0~,
  \label{eq.semigroup.calL_3}
\end{equation}
see \cite[Appendix~A]{GW2}. Note that the semigroup $e^{t\LL_3}$ is 
not strongly continuous in the space $BC(\R)$ equipped with the 
supremum norm. This is mainly due to the dilation factor $e^{-t}$ in 
\eqref{eq.semigroup.calL_3}. However, if we equip $BC(\R)$ with 
the (weaker) topology of uniform convergence on compact sets, then 
the map $t \mapsto e^{t\LL_3}f$ is continuous for any $f \in BC(\R)$. 
This observation is the reason for introducing the space $X_{loc}(m)$
in Section~\ref{sec.intro}.

Since the operators $\LL_h$ and $\LL_3$ act on different variables, 
it is easy to obtain the semigroup generated by $L_3 = \LL_h + 
\LL_3$ by combining the formulas \eqref{eq.semigroup.calL_h} 
and \eqref{eq.semigroup.calL_3}. We find 
\begin{equation}
  (e^{tL_3}\phi)(x) \,=\, \frac{1}{\sqrt{2\pi a(2t)}}\int_{\R}
  e^{-\frac{|x_3e^{-t}-y_3|^2}{2a(2t)}}\Bigl(e^{t\LL_h}\phi(\cdot,y_3)
  \Bigr)(x_h)\dd y_3~, \qquad t > 0~. \label{eq.semigroup.L_3}
\end{equation}
In \cite[Proposition A.6]{GW2}, it is shown that this expression
defines a uniformly bounded semigroup in $X(m)$ for any $m > 1$, and
that the map $t \mapsto e^{tL_3}$ is strongly continous in the
topology of $X_{loc}(m)$. Moreover, the subspace $X_0(m)$ is left
invariant by $e^{tL_3}$ for any $t \ge 0$.  Using these results and
the relation \eqref{eq.Lexp}, we conclude that the
three-dimensional operator $L$ generates a uniformly bounded semigroup
in the space $\X(m)$, given by
\begin{equation}
  e^{tL}\omega \,=\, \Bigl(e^{-3t/2}e^{tL_3}\omega_1\,,\,
  e^{-3t/2}e^{tL_3}\omega_2\,,\,e^{tL_3}\omega_3\Bigr)^\top~, 
  \qquad t \ge 0~. \label{eq.semigroup.L}
\end{equation}
As is easily verified, if $\nabla \cdot \omega = 0$, then 
$\nabla\cdot e^{tL}\omega = 0$ for all $t \ge 0$. 

The asymptotic stability of the Burgers vortices relies heavily
on the decay properties of the semigroup $e^{tL}$ as $t \to \infty$. 
In the proof of Theorems~\ref{thm.main1} and \ref{thm.main2}, 
we also use the smoothing properties of the operator $e^{tL}$ for 
$t > 0$, and in particular the fact that $e^{tL}$ extends to 
a bounded operator from $\X^p(m)$ into $\X^2(m)$ for all $p \in [1,2]$. 
All the needed estimated are collected in the following statement. 

\begin{prop}\label{prop.semigroup.L} Let $m\in (1,\infty]$, 
 $p\in [1,2]$, and take $\eta \in (0,\frac12]$ such that $2\eta < 
m - 1$. For any $\beta = (\beta_1,\beta_2,\beta_3) \in \N^3$, there 
exists $C > 0$ such that the following estimates hold:
\begin{align}
  \|\partial_x^\beta e^{tL_h} \omega_h \|_{X(m)^2} \,&\le\, \frac{C
  e^{-(\frac{3}{2}+\beta_3) t}}{a(t)^{\frac{1}{p}-\frac{1}{2}+\frac{|\beta|}{2}}}
  \|\omega_h\|_{X^p(m)^2}~, \label{est.prop.semigroup.L.1}\\
  \|\partial_x^\beta e^{tL_3} \omega_3\|_{X(m)} \,&\le\, \frac{C
  e^{-(\eta+\beta_3)t}}{a(t)^{\frac{1}{p}-\frac{1}{2}+\frac{|\beta|}{2}}}
  \|\omega_3\|_{X^p(m)}~, \label{est.prop.semigroup.L.2}
\end{align}
for any $\omega \in \X^p(m)$ and all $t > 0$. Here $a(t) = 1-e^{-t}$ 
and $|\beta| = \beta_1 + \beta_2 + \beta_3$. 
\end{prop}

\noindent{\bf Proof.} We first assume that $m \in (1,\infty)$. If 
$p \in [1,2]$ and $\beta_h = (\beta_1,\beta_2) \in \N^2$, it is 
proved in \cite[Appendix~A]{GW0} that 
\begin{equation}
  \|\partial_{x_h}^{\beta_h} e^{t\LL_h} f\|_{L^2(m)} \,\le\,
  \frac{C}{a(t)^{\frac{1}{p}-\frac{1}{2}+\frac{|\beta_h|}{2}}}
  \|f\|_{L^p(m)}~, \qquad t > 0~,\label{est.semigroup_h.1}
\end{equation}
for all $f \in L^p(m)$. If in addition $f \in L^p_0(m)$, we have the 
stronger estimate 
\begin{equation}
  \|\partial_{x_h}^{\beta_h} e^{t\LL_h} f\|_{L^2(m)} \,\le\,
  \frac{C e^{-\eta t}}{a(t)^{\frac{1}{p}-\frac{1}{2}+\frac{|\beta_h|}{2}}}
  \|f\|_{L^p(m)}~, \qquad t > 0~,\label{est.semigroup_h.2}
\end{equation}
where $\eta > 0$ is as in Proposition~\ref{prop.semigroup.L}. On the 
other hand, using \eqref{eq.semigroup.L_3}, we find by direct
calculation
\begin{equation}
  \|\partial_{x_3}^{\beta_3}e^{t\LL_3}f\|_{L^\infty(\R)} \,\le\,
  \frac{C e^{-\beta_3 t}}{a(t)^\frac{\beta_3}{2}}
  \|f\|_{L^\infty(\R)}~, \qquad t > 0~. \label{est.semigroup_3}
\end{equation}
Here, as in \eqref{eq.semigroup}, the stabilizing factor
$e^{-\beta_3t}$ comes from the dilation operator $-x_3\partial_{x_3}$
which enters the definition of $\LL_3$. Now, if we start from the
representation \eqref{eq.semigroup.L_3} and use the estimates
\eqref{est.semigroup_h.1}--\eqref{est.semigroup_3}, we easily 
obtain \eqref{est.prop.semigroup.L.1}, \eqref{est.prop.semigroup.L.2}
by a direct calculation, see \cite[Proposition A.6]{GW2}. 

To complete the proof of Proposition~\ref{prop.semigroup.L}, it
remains to show that \eqref{est.semigroup_h.1},
\eqref{est.semigroup_h.2} still hold when $m = \infty$. If $t \in
(0,1)$, estimate \eqref{est.semigroup_h.1} is easily obtained by a
direct calculation, based on the representation
\eqref{eq.semigroup.calL_h}. Using this remark and the semigroup
property of $e^{t\LL_h}$, we conclude that it is sufficient to
establish \eqref{est.semigroup_h.1}, \eqref{est.semigroup_h.2} in the
particular case where $p = 2$ and $\beta_h = 0$. This in turns follows
easily from the spectral properties of the generator $\LL_h$. Indeed,
it is well-known that $\LL_h$ is a self-adjoint operator in
$L^2(\infty)$ with purely discrete spectrum $\sigma(\LL_h) = 
\{-\frac{k}{2}\,|\,k=0,1,2,\dots\}$. Moreover, the subspace 
$L^2_0(\infty)$ is precisely the orthogonal complement of the 
eigenspace corresponding to the zero eigenvalue, see for 
example \cite[Lemma 4.7]{GW1}. It follows that $e^{t\LL_h}$ 
is a semigroup of contractions in $L^2(\infty)$, and that 
$\|e^{t\LL_h}f\|_{L^2(\infty)} \le e^{-t/2}\|f\|_{L^2(\infty)}$ 
for all $t \ge 0$ if $f \in L^2_0(\infty)$. This proves  
\eqref{est.semigroup_h.1} and \eqref{est.semigroup_h.2}, with 
$\eta = 1/2$. \QED

\subsection{Estimates for the velocity fields}\label{subsec.est.velocity}

If the velocity $u$ and the vorticity $\omega$ are related by the 
Biot-Savart law \eqref{def.3DBS}, we have $|u| \le J(|\omega|)$, where 
$J$ is the Riesz potential defined by
\begin{equation}
  J(\phi)(x) \,=\, \frac{1}{4\pi}\int_{\R^3}\frac{1}{|x-y|^2}
  \,\phi(y)\dd y~, \qquad x \in \R^3~. \label{def.Riesz}
\end{equation}
Since $\omega$ will typically belong to the Banach space 
$\X(m)$, we need estimates on the Riesz potential $J(\phi)$ 
for $\phi \in X(m)$. We start with a preliminary result:

\begin{lem}\label{lem.Riesz} Let $p_1 \in [1,2)$, $p_2 \in [1,2]$, 
and assume that $\phi\in X^{p_1}(0)\cap X^{p_2}(0)$. If 
$q_1, q_2 \in [1,\infty]$ satisfy
\begin{equation}\label{def.q1q2}
  \frac{2p_1}{2-p_1} \,<\, q_1 \,\le\, \infty~, \qquad
  p_2 \,<\, q_2 \,<\, \frac{2p_2}{2-p_2}~, 
\end{equation}
then $J(\phi) = J_1(\phi)+J_2(\phi)$ with $J_i(\phi)\in X^{q_i}(0)$ 
for $i = 1,2$, and we have the following estimates
\begin{align}
  \|J_1(\phi)\|_{X^{q_1}(0)}  \,&\le\, C(p_1,q_1)\|\phi\|_{X^{p_1}(0)}
  ~,\label{est.J1}\\
  \|J_2(\phi)\|_{X^{q_2}(0)}  \,&\le\, C(p_2,q_2)\|\phi\|_{X^{p_2}(0)}
  ~.\label{est.J2}
\end{align}
\end{lem}

\noindent{\bf Proof.} We proceed as in \cite[Proposition A.9]{GW2}. 
We first observe that
\begin{align*}
  J(\phi)(x_h,x_3) \,&=\, \int_{|x_3-y_3|\ge 1}F(x_h;x_3,y_3)\dd y_3 + 
  \int_{|x_3-y_3|< 1}F(x_h;x_3,y_3)\dd y_3 \\ 
  \,&=\, J_1(\phi)(x_h,x_3) + J_2(\phi)(x_h,x_3)~,
\end{align*}
where
\[
  F(x_h;x_3,y_3) \,=\, \int_{\R^2}\frac{\phi(y_h,y_3)}{|x_h-y_h|^2
  +(x_3-y_3)^2}\dd y_h~, \qquad x_h \in \R^2~, \quad x_3,y_3 \in \R~. 
\]
For any $a \in \R$, let $f_a(y_h) = (a^2+|y_h|^2)^{-1}$. Then 
$f_a \in L^r(\R^2)$ for any $r > 1$ and any $a \neq 0$, and 
there exists $C_r > 0$ such that
\[
  \|f_a\|_{L^r(\R^2)} \,\le\, \frac{C_r}{|a|^{2-\frac2r}}~.
\]
Moreover, we have $F(\cdot\,;x_3,y_3) = \phi(\cdot,y_3)\star f_{x_3-y_3}$
by construction. Thus, if we take $1 \le p,q,r \le \infty$ such that
$1+\frac1q = \frac1p+\frac1r$, we obtain using Young's inequality
\[
  \|F(\cdot\,;x_3,y_3)\|_{L^q(\R^2)} \,\le\, 
  \|\phi(\cdot,y_3)\|_{L^p(\R^2)} \|f_{x_3-y_3}\|_{L^r(\R^2)}
  \,\le\, \frac{C_r \|\phi(\cdot,y_3)\||_{L^p(\R^2)}}{|x_3-
  y_3|^{2-\frac2r}}~.
\]

To estimate $J_1(\phi)$, we choose $p = p_1$, $q = q_1$. 
In view of \eqref{def.q1q2}, the corresponding exponent
$r = r_1$ satisfies $2 < r_1 \le \infty$, so that $2-\frac2{r_1}
\in (1,2]$. By Minkowski's inequality, we thus find
\[
  \|J_1(\phi)(\cdot,x_3)\|_{L^{q_1}(\R^2)} \,\le\, \int_{|x_3-y_3|\ge 1}
  \|F(\cdot\,;x_3,y_3)\|_{L^{q_1}(\R^2)}\dd y_3 \,\le\, 
  C(r_1) \sup_{y_3 \in \R}\|\phi(\cdot,y_3)\|_{L^{p_1}(\R^2)}~.
\]
Taking the supremum over $x_3 \in \R$, we obtain \eqref{est.J1}. 
Similarly, to bound $J_2(\phi)$, we take  $p = p_2$, $q = q_2$. 
Then $1 < r_2 < 2$, so that $2 - \frac2{r_2} \in (0,1)$. 
We thus obtain
\[
  \|J_2(\phi)(\cdot,x_3)\|_{L^{q_2}(\R^2)} \,\le\, \int_{|x_3-y_3|< 1}
  \|F(\cdot\,;x_3,y_3)\|_{L^{q_2}(\R^2)}\dd y_3 \,\le\, 
  C(r_2) \sup_{y_3 \in \R}\|\phi(\cdot,y_3)\|_{L^{p_2}(\R^2)}~, 
\]
and \eqref{est.J2} follows. Finally, the uniform continuity of 
$J_i(\phi)(\cdot,x_3)$ with respect to $x_3$ can be verified exactly as 
in the proof of \cite[Proposition~A.9]{GW2}. \QED

\medskip

As an immediate consequence, we obtain the following useful statements.

\begin{prop}\label{prop.Riesz} Let $\phi \in X(m)$ for some $m \in 
(1,\infty]$. Then $J(\phi)\in X^q(0)$ for all $q\in (2,\infty)$, and 
there exists a positive constant $C = C(m,q)$ such that
\begin{equation}
  \|J(\phi)\|_{X^q(0)} \,\le\, C\|\phi\|_{X(m)}~.
  \label{est.prop.Riesz}
\end{equation}
\end{prop} 

\noindent{\bf Proof.} If $m > 1$, we recall that $X(m) \hookrightarrow
X^p(0)$ for all $p \in [1,2]$. Thus we can apply Lemma~\ref{lem.Riesz} 
with $p_1 = 1$, $p_2 = 2$, and $q_1 = q_2 = q \in (2,\infty)$, 
and the result follows. \QED

\begin{cor}\label{cor.prop.Riesz} Let $\phi_1$, $\phi_2 \in X(m)$ for 
some $m \in (1,\infty]$. Then $\phi_1J(\phi_2)\in X^p(m)$ for all 
$p\in (1,2)$, and there exists a positive constant $C = C(m,p)$ such that 
\begin{equation}
  \|\phi_1J(\phi_2)\|_{X^p(m)} \,\le\, C\|\phi_1\|_{X(m)}
  \|\phi_2\|_{X(m)}~.
  \label{est.cor.prop.Riesz}
\end{equation}
\end{cor} 

\noindent{\bf Proof.} We proceed as in \cite[Corollary~A.10]{GW2}. 
Let $p \in (1,2)$, and take $q \in (2,\infty)$ such that 
$\frac1q = \frac1p -\frac12$. For any $x_3 \in \R$, we have by 
H\"older's inequality
\begin{align*}
  \|\phi_1(\cdot,x_3)J(\phi_2)&(\cdot,x_3)\|_{L^p(m)} \,=\, 
  \Bigl(\int_{\R^2}\rho_m(|x_h|^2)^{p/2}|\phi_1(x_h,x_3)|^p |J(\phi_2)
  (x_h,x_3)|^p\dd x_h\Bigr)^{1/p} \\
  \,&\le\, \Bigl(\int_{\R^2}\rho_m(|x_h|^2)|\phi_1(x_h,x_3)|^2\dd x_h
  \Bigr)^{1/2}\Bigl(\int_{\R^2}|J(\phi_2)(x_h,x_3)|^q\dd x_h\Bigr)^{1/q} \\
  \,&=\, \|\phi_1(\cdot,x_3)\|_{L^2(m)} \|J(\phi_2)(\cdot,x_3)\|_{L^q(0)}~.
\end{align*}
Taking the supremum over $x_3 \in \R$ and using
Proposition~\ref{prop.Riesz}, we obtain \eqref{est.cor.prop.Riesz}.
Finally, it is clear that the map $x_3 \mapsto \phi_1(\cdot,x_3)
J(\phi_2)(\cdot,x_3)$ is continuous from $\R$ into $L^p(m)$. \QED

\medskip

We conclude this section with an estimate on the linear operator
\eqref{def.Lambda} which will be needed in Section~\ref{sec.linear}.

\begin{lem}\label{lem.after.Riesz} Let $p \in [1,2]$ and $2-\frac2p < m
\le \infty$. For any $\beta \in \N^3$, there exists $C > 0$ such that
\begin{equation}
  \|\partial_x^\beta \Lambda \omega\|_{\X^p (m)} ~\le~ 
  C \!\!\sum_{|\tilde\beta|\le |\beta|+1} \|\partial_x^{\tilde \beta} 
  \omega\|_{\X^p(m)}. \label{est.lem.after.Riesz}
\end{equation}
\end{lem}

\noindent{\bf Proof.} It is sufficient to prove 
\eqref{est.lem.after.Riesz} for $\beta = 0$. The general 
case easily follows if we use the Leibniz rule to differentiate 
$\Lambda \omega$ (we omit the details). 

Assume thus that $\omega$ belongs to $\X^p(m)$, together with
its first order derivatives. Since the function $U^G$ defined
in \eqref{def.u^g} is smooth and bounded (together with 
all its derivatives), it is clear that
\[
  \|(U^G,\nabla)\omega\|_{\X^p(m)} + \|(\omega,\nabla)U^G\|_{\X^p(m)}
  \,\le\, C\sum_{|\tilde \beta|\le 1}\|\partial_x^{\tilde \beta}\omega
  \|_{\X^p(m)}~.
\]
We now estimate the term $(u,\nabla)G = (K_{3D}*\omega,\nabla) G$, 
using the fact that $|K_{3D}*\omega| \le J(|\omega|)$. Since 
$|\omega| \in X^1(0) \cap X^p(0)$ by assumption, we can apply
Lemma~\ref{lem.Riesz} with $p_1 = 1$, $q_1 = \infty$, $p_2 = p$, 
and $q_2\in (p,\frac{2p}{2-p})$. By H\"older's inequality,
we easily find
\begin{align*}
  \|J_1(|\omega|) |\nabla G|\|_{X^p(m)} \,&\le\, C\|J_1(|\omega|)
  \|_{X^\infty(0)} \,\le\, C\||\omega|\|_{X^1(0)} \,\le\, 
  C\|\omega\|_{\X^p(m)}~,\\
  \|J_2(|\omega|) |\nabla G|\|_{X^p(m)} \,&\le\, C\|J_2(|\omega|)
  \|_{X^{q_2}(0)} \,\le\, C\||\omega|\|_{X^p(0)} \,\le\, 
  C\|\omega\|_{\X^p(m)}~.
\end{align*}
We conclude that $\|(u,\nabla)G\|_{\X^p(m)} \le \|(K_{3D}*\omega,
\nabla)G\|_{\X^p(m)} \le C\|\omega\|_{\X^p(m)}$. In a similar way, 
commuting the derivative and the convolution operator, we obtain 
the estimate $\|(G,\nabla)u\|_{\X^p(m)} \le \|(G,\nabla)(K_{3D}
*\omega)\|_{\X^p(m)} \le C\|\nabla \omega\|_{\X^p(m)}$. This completes 
the proof. \QED


\section{The vectorial 2D problem}\label{sec.vectorial}

In this section we study the linearized equation $\partial_t \omega 
= (L-\alpha\Lambda)\omega$ in the particular case where the vorticity
$\omega$ does not depend on the vertical variable. As was explained
in the introduction, this preliminary step is an essential ingredient
in the linear stability proof which will be presented in 
Section~\ref{sec.linear}. 

If $\partial_{x_3}\omega = 0$, then $\LL_3\omega = 0$, and the 
expression \eqref{eq.Lexp} of the linear operator $L$ becomes 
significantly simpler. On the other hand, we know from 
\eqref{def.Lambda} that
\begin{equation}
  \Lambda \omega \,=\, \Lambda_1\omega - \Lambda_2\omega + 
  \Lambda_3\omega-\Lambda_4\omega~,
  \label{eq.Lambda_dec}
\end{equation}
where
\begin{equation}
  \begin{array}{l}
  \Lambda_1 \omega \,=\, (U^G,\nabla)\omega \,=\, (U^G_h,\nabla_h)
  \omega~, \\[2mm]
  \Lambda_2 \omega \,=\, (\omega,\nabla)U^G \,=\, 
  (\omega_h,\nabla_h)U^G~,  
  \end{array} \qquad
  \begin{array}{l}
  \Lambda_3 \omega \,=\, (u,\nabla)G \,=\, (u_h,\nabla_h)G~, \\[2mm]
  \Lambda_4 \omega \,=\, (G,\nabla)u \,=\, g\partial_{x_3}u~.  
  \end{array}
  \label{def.Lambda_i}
\end{equation}
Here $u = K_{3D}*\omega$ is the velocity field obtained from $\omega$
via the three-dimensional Biot-Savart law \eqref{def.3DBS}. Since 
$\partial_{x_3}\omega = 0$, we have $\partial_{x_3}u = 0$, hence 
$\Lambda_4 \omega = 0$ in our case. Moreover, it is easy to verify
that $u = (u_h,u_3)$, where $u_h = K_{2D}\star\omega_3$. Thus, we see 
that
\begin{equation}
  (L-\alpha\Lambda)\omega \,=\, \cL_\alpha \omega \,:=\,
  \begin{pmatrix} (\LL_h-\frac32)\omega_h - \alpha (\Lambda_1 - 
  \tilde\Lambda_2)\omega_h\\
  \LL_h \omega_3 - \alpha (\Lambda_1 + \tilde\Lambda_3)\omega_3
  \end{pmatrix} \,\equiv\, 
  \begin{pmatrix} \cL_{\alpha,h}\,\omega_h \\ \cL_{\alpha,3}\,\omega_3 
  \end{pmatrix}~,
  \label{def.calLambda}
\end{equation}
where $\tilde \Lambda_2 \omega_h = (\omega_h,\nabla_h)U_h^G$ and
$\tilde \Lambda_3 \omega_3 = (K_{2D}\star\omega_3,\nabla_h)g$.  

For any $\alpha \in \R$ and any $m \in (1,\infty]$, the operator
$\cL_\alpha$ defined by \eqref{def.calLambda} is the generator of a
strongly continuous semigroup in the space $L^2(m)^3$. This property
can be established by a standard perturbation argument, see
Lemma~\ref{prop.vectorial.2D.local} below. Our main goal here is to
obtain accurate decay estimates for the semigroup $e^{t\cL_\alpha}$ as
$t \to \infty$. As is clear from \eqref{def.calLambda}, the evolutions
for $\omega_h$ and $\omega_3$ are completely decoupled, so that we
can consider the semigroups $e^{t\cL_{\alpha,h}}$ and $e^{t\cL_{\alpha,3}}$ 
separately. The main contribution of this section is:

\begin{prop}\label{prop.vectorial.2D} Fix $m \in (1,\infty]$, 
$\alpha \in \R$, $\mu \in (0,\frac32)$, and take $\eta \in (0,\frac12]$ 
such that $1+2\eta < m$. Then there exists $C > 0$ such that 
\begin{align}
  \|e^{t\cL_{\alpha,h}}\omega_h\|_{L^2(m)^2} \,&\le\, C\,e^{-\mu t}
  \|\omega_h\|_{L^2(m)^2}~, \qquad t \ge 0~,
  \label{est.prop.vectorial.2D.1} \\
  \|e^{t\cL_{\alpha,3}}\omega_3\|_{L^2(m)} \,&\le\, C\,e^{-\eta t}
  \|\omega_3\|_{L^2(m)}~, ~\qquad t \ge 0~,
  \label{est.prop.vectorial.2D.2}
\end{align}
for all $\omega \in L^2(m)^2 \times L^2_0(m)$. 
\end{prop}

Estimate \eqref{est.prop.vectorial.2D.2} was obtained in
\cite[Proposition~4.12]{GW1} for $m < \infty$, and the proof given
there extends to the limiting case $m = \infty$ without additional
difficulty. Remark that the decay rate $e^{-\eta t}$ is obtained 
using the fact that $\omega_3 \in L^2_0(m)$: If we only assume that 
$\omega_3 \in L^2(m)$ for some $m > 1$, then 
\eqref{est.prop.vectorial.2D.2} holds with $\eta = 0$. 
Note, however, that $\omega$ is not assumed to be divergence-free
in this section.

From now on, we focus on the semigroup $e^{t\cL_{\alpha,h}}$, which
has not been studied yet. To prove \eqref{est.prop.vectorial.2D.1}, 
we use the same arguments as in \cite[Section~4.2]{GW1}. We 
first establish a short time estimate: 

\begin{lem}\label{prop.vectorial.2D.local} Fix $m \in (1,\infty]$, 
$\alpha \in \R$, and $T > 0$. There exists $C = C(T,m,|\alpha|) > 0$ 
such that  
\begin{equation}
  \sup_{0 \le t \le T}\Bigl(\|e^{t\cL_{\alpha,h}}\omega_h\|_{L^2(m)^2}
  +a(t)^\frac{1}{2}\|\nabla_h e^{t\cL_{\alpha,h}}\omega_h\|_{L^2(m)^4}
  \Bigr) \,\le\, C\|\omega_h\|_{L^2(m)^2}~,
  \label{est.lem.local}
\end{equation}
for all $\omega_h \in L^2(m)^2$. Here $a(t) = 1-e^{-t}$. 
\end{lem}

\noindent{\bf Proof.} Given $\omega_h^0 \in L^2(m)^2$, the idea
is to solve the integral equation
\begin{equation}
  \omega_h(t) \,=\, e^{t(\LL_h-\frac{3}{2})}\omega_h^0 -\alpha
  \int_0^t e^{(t-s)(\LL_h-\frac{3}{2})}(\Lambda_1-\tilde\Lambda_2) 
  \omega_h(s)\dd s~, \qquad t \in [0,T]~,
  \label{eq.integ.2Dvec}
\end{equation}
by a fixed point argument in the space $X_T = \{\omega_h \in C([0,T],
L^2(m)^2\,|\, \|\omega_h\|_{X_T} < \infty\}$ defined by the norm
\[
  \|\omega_h\|_{X_T} \,=\, \sup_{0\le t\le T}\|\omega_h(t)\|_{L^2(m)^2}
  + \sup_{0\le t\le T}a(t)^\frac{1}{2}\|\nabla_h \omega_h(t)
  \|_{L^2(m)^4}~.
\]
From \eqref{est.semigroup_h.1} we know that $\|e^{t(\LL_h-\frac{3}{2})}
\omega_h^0\|_{X_T} \le C_1\|\omega_h^0\|_{L^2(m)^2}$, for some $C_1 > 0$ 
independent of $T$. To estimate the integral term in 
\eqref{eq.integ.2Dvec}, we first observe that the velocity field 
$U^G$ defined by \eqref{def.u^g} satisfies
\begin{equation}\label{eq.UGbounds}
  \sup_{x_h \in \R^2}(1+|x_h|)|U^G(x_h)| + \sup_{x_h \in \R^2}(1+|x_h|)^2
  |\nabla_h U^G(x_h)| \,<\, \infty~.
\end{equation}
In view of the definitions \eqref{def.Lambda_i}, we thus have
\begin{align}
   \|(1+|x_h|) \Lambda_1 \omega_h\|_{L^2(m)^2} \,&\le\, 
   C\|\nabla_h \omega_h\|_{L^2(m)^4}~, 
   \label{est.Lambda_1.subsec.2D}\\
   \|(1+|x_h|)^2 \tilde\Lambda_2 \omega_h\|_{L^2(m)^2}\,&\le\, 
   C\|\omega_h\|_{L^2(m)^2}~.
   \label{est.Lambda_2.subsec.2D}
\end{align}
Using these estimates together with \eqref{est.semigroup_h.1}, 
we can bound
\begin{align*}
  \Bigl\|\int_0^t e^{(t-s)(\LL_h-\frac{3}{2})}&(\Lambda_1-
  \tilde \Lambda_2)\omega_h(s) \dd s\Bigr\|_{L^2(m)^2} \\ 
  \,&\le\, C \int_0^t e^{-\frac32(t-s)}\Bigl(\|\omega_h(s)\|_{L^2(m)^2} + 
  \|\nabla_h \omega_h(s)\|_{L^2(m)^4}\Bigr)\dd s \\
  \,&\le\, C \|\omega_h\|_{X_T} \int_0^t e^{-\frac32(t-s)} 
  a(s)^{-\frac12}\dd s \,\le\, C a(T)^{\frac12}\|\omega_h\|_{X_T}~. 
\end{align*}
In a similar way, 
\begin{align}
  \Bigl\|\nabla_h\int_0^t &e^{(t-s)(\LL_h-\frac{3}{2})}(\Lambda_1-\tilde
  \Lambda_2) \omega_h(s)\dd s\Bigr\|_{L^2(m)^4} \label{est.nabla_int} \\
  \,&\le\, C \int_0^t \frac{e^{-\frac32(t-s)}}{a(t-s)^{\frac12}} 
  \Bigl(\|\omega_h(s)\|_{L^2(m)^2} + \|\nabla_h \omega_h(s)\|_{L^2(m)^4}
  \Bigr)\dd s \,\le\, C \|\omega_h\|_{X_T}~. 
  \nonumber 
\end{align}
Summarizing, we have shown that $\|\omega_h\|_{X_T} \le C_1
\|\omega_h^0\|_{L^2(m)^2} + C_2 |\alpha| a(T)^{1/2}\|\omega_h\|_{X_T}$, 
for some positive constants $C_1, C_2$. If we now take $T > 0$ small
enough so that $C_2 |\alpha| a(T)^{1/2} \le 1/2$, we see that the 
right-hand side of \eqref{eq.integ.2Dvec} is a strict contraction
in $X_T$. We deduce that \eqref{eq.integ.2Dvec} has a unique solution, 
which satisfies $\|\omega_h\|_{X_T} \le 2C_1\|\omega_h^0\|_{L^2(m)^2}$. 
Since $\omega_h(t) = e^{t\cL_{\alpha,h}}\omega_h^0$ by construction, 
this proves \eqref{est.lem.local} for $T$ sufficiently small, 
and the general case follows due to the semigroup property. 
This concludes the proof. \QED

\medskip We next consider the essential spectrum of the semigroup
$e^{t\cL_{\alpha,h}}$, and begin with a few definitions.  If $A$ is a
bounded linear operator on a (complex) Banach space $X$, we define the
essential spectrum $\sigma_{ess}(A\,;X)$ as the set of all $z \in \C$
such that $A-z$ is not a Fredholm operator with zero index, see
\cite{Ka} or \cite{EN}. The essential spectral radius of $A$ in $X$ is
given by
\[
  r_{ess}(A\,;X) \,=\, \sup\Bigl\{|z|\,;\, z \in \sigma_{ess}(A\,;X)
  \Bigr\} \,<\,\infty~.
\]
If $|z| > r_{ess}(A\,;X)$, then either $z$ is in the resolvent set of 
$A$, or $z$ is an eigenvalue of $A$ with finite multiplicity, see 
\cite[Corollary~IV.2.11]{EN}. In the latter case, we say that $z$ 
belongs to the discrete spectrum of $A$. 

In what follows, we consider the linear operator $\cL_{\alpha,h}$ 
as acting on the complexified space $L^2(m)^2$, i.e. the space 
of all $\omega_h : \R^2 \to \C^2$ such that $\|\omega_h\|_{L^2(m)^2}
< \infty$. Our first result shows that the essential spectral 
radius of the operator $e^{t\cL_{\alpha,h}}$ in $L^2(m)^2$ does not 
depend on $\alpha$. 

\begin{prop}\label{prop.essential.spectrum} Let $m\in (1,\infty]$ and 
$\alpha \in \R$. Then for each $t > 0$ we have
\begin{equation}
  r_{ess}\Bigl(e^{t\cL_{\alpha,h}}\,;L^2(m)^2\Bigr) \,=\, 
  r_{ess}\Bigl(e^{t\cL_{0,h}}\,;L^2(m)^2\Bigr) \,=\, 
  e^{-(\frac{m}{2}+1)t}~.
  \label{est.prop.essential.spectrum.1}
\end{equation}
\end{prop}

\noindent{\bf Proof.} Since $\cL_{0,h} = \LL_h-\frac32$, the last
equality in \eqref{est.prop.essential.spectrum.1} follows from
\cite[Theorem A.1]{GW0} if $m < \infty$. If $m = \infty$, then
$e^{t\LL_h}$ is a compact operator for any $t > 0$, hence
$r_{ess}(e^{t\cL_{0,h}}\,;L^2(\infty)^2) = 0$. To prove the first
equality in \eqref{est.prop.essential.spectrum.1}, we fix $t > 0$. Our
goal is to show that the linear operator $\Delta_\alpha(t) =
e^{t\cL_\alpha,h} - e^{t(\LL_h-\frac32)}$ is compact in $L^2(m)^2$. By
Weyl's theorem, this will imply that both semigroups have the same
essential spectrum, hence the same essential spectral radius.  In view
of \eqref{eq.integ.2Dvec} we have, for all $\omega_h \in L^2(m)^2$,
\begin{equation}
  \Delta_\alpha(t)\omega_h \,=\, -\alpha \int_0^t e^{(t-s)(\LL_h
  -\frac{3}{2})}(\Lambda_1-\tilde\Lambda_2)e^{s\cL_{\alpha,h}} 
  \omega_h\dd s~.
  \label{def.Delta}
\end{equation}
Let $w(x_h) = 1+|x_h|$. If $m < \infty$, it follows from 
\eqref{est.semigroup_h.1} and definition \eqref{def.L^2(m)} 
that
\begin{equation}
  \|w\,e^{t\LL_h}\omega_h\|_{L^2(m)^2} \,\le\, 
  C \|e^{t\LL_h}\omega_h\|_{L^2(m+1)^2} \,\le\, 
  C\|w\,\omega_h\|_{L^2(m)^2}~, 
  \label{eq.weightLL_h.1}
\end{equation}
for all $\omega_h \in L^2(m)^2$ and all $t \ge 0$. If $m = \infty$, we
know from \cite[Proposition~2.1]{GW3} that $w (-\LL_h+1)^{-1/2}$ is a
bounded operator in $L^2(\infty)^2$, and since $\LL_h$ is the
generator of an analytic semigroup we easily obtain
\begin{equation}
  \|w\,e^{t\LL_h}\omega_h\|_{L^2(m)^2} \,\le\, 
  C \|(-\LL_h+1)^{1/2}e^{t\LL_h}\omega_h\|_{L^2(m)^2} \,\le\, 
  \frac{C}{a(t)^{1/2}}\,\|\omega_h\|_{L^2(m)^2}~,
  \label{eq.weightLL_h.2}
\end{equation}
for all $t > 0$. Now, starting from \eqref{def.Delta} and 
using either \eqref{eq.weightLL_h.1} or \eqref{eq.weightLL_h.2} 
together with \eqref{est.Lambda_1.subsec.2D}, 
\eqref{est.Lambda_2.subsec.2D}, and Lemma~\ref{prop.vectorial.2D.local}, 
we find
\begin{align*}
  \|w\,\Delta_\alpha(t)\omega_h\|_{L^2(m)^2} \,&\le\, C|\alpha|\int_0^t 
  \frac{e^{-\frac32(t-s)}}{a(t{-}s)^{1/2}}\Bigl(\|e^{s\cL_{\alpha,h}}
  \omega_h\|_{L^2(m)^2} + \|\nabla_h e^{s\cL_{\alpha,h}}\omega_h
  \|_{L^2(m)^4}\Bigr)\dd s \\
  \,&\le\, C |\alpha| \|\omega_h\|_{L^2(m)^2} \int_0^t \frac{
  e^{-\frac32(t-s)}}{a(t{-}s)^{1/2}a(s)^{1/2}}\dd s \,\le\, 
  C |\alpha| \|\omega_h\|_{L^2(m)^2}~. 
\end{align*}
Moreover, proceeding as in \eqref{est.nabla_int}, we find
$\|\nabla_h \Delta_\alpha(t) \omega_h\|_{L^2(m)^4} \le C |\alpha| 
\|\omega_h \|_{L^2(m)^2}$. Thus we have shown that $w \Delta_\alpha(t)$ 
and $\nabla_h \Delta_\alpha(t)$ are bounded operators in $L^2(m)$. By
Rellich's criterion, we conclude that $\Delta_\alpha(t)$ is a compact 
operator in $L^2(m)^2$, for any $t > 0$. This completes the proof.  
\QED

\medskip In view of Proposition~\ref{prop.essential.spectrum}, the
spectrum of the semigroup $e^{t\cL_{\alpha,h}}$ outside the disk of
radius $e^{-(\frac{m}{2}+1)t}$ in the complex plane is purely
discrete. By the spectral mapping theorem \cite{EN}, to control that
part of the spectrum it is sufficient to locate the eigenvalues of the
generator $\cL_{\alpha,h}$. Thus we look for nontrivial solutions of
the eigenvalue problem
\begin{equation}
  \cL_{\alpha,h}\omega_h \,=\, \lambda \omega_h~,
  \label{eq.eigenvalue}
\end{equation}
where $\omega_h \in L^2(m)^2$ and $\lambda \in \C$ satisfies 
$\Re\lambda > -\frac{m}{2} - 1$. The following auxiliary result 
shows that the eigenfunctions $\omega_h$ always have a Gaussian 
decay at infinity. 

\begin{prop}\label{prop.discrete.spectrum} 
Let $m \in (1,\infty)$ and $\alpha\in\R$. If $\omega_h \in L^2(m)^2$ 
is a solution of \eqref{eq.eigenvalue} with $\Re\lambda > -\frac{m}{2}
-1$, then $\omega_h\in L^2(\infty)^2$. 
\end{prop}

The proof of Proposition~\ref{prop.discrete.spectrum} is postponed
to Section~\ref{subsec.discrete} below. Note that a similar result 
for the nonlocal operator $\cL_{\alpha,3}$ has been obtained in 
\cite[Lemma 4.5]{GW1}, and plays a key role in the derivation
of estimate~\eqref{est.prop.vectorial.2D.2}. Thanks to 
Proposition~\ref{prop.discrete.spectrum}, we only need to control 
the eigenvalues of $\cL_{\alpha,h}$ in the Gaussian space $L^2(\infty)^2$. 
This is the last important step in the proof of 
Proposition~\ref{prop.vectorial.2D}. 

\begin{prop}\label{prop.discrete.spectrum.2} If $\lambda$ is an 
eigenvalue of $\cL_{\alpha,h}$ in $L^2(\infty)^2$, then $\Re\lambda 
\le -\frac32$.
\end{prop}

\noindent{\bf Proof.} Assume that $\omega_h \in L^2(\infty)^2$ 
is a nontrivial solution of the eigenvalue problem \eqref{eq.eigenvalue}, 
for some $\alpha \in \R$ and some $\lambda \in \C$. Using 
\eqref{def.calLambda}, we thus have
\begin{equation}\label{eq.eigenvalue2}
  \lambda \omega_h \,=\, \LL_h \omega_h -\frac32 \omega_h 
  -\alpha (U^G_h,\nabla_h)\omega_h + \alpha (\omega_h,\nabla_h)U^G_h~,
\end{equation}
where the velocity field $U^G$ is defined in \eqref{def.u^g}. Since
$\cL_{\alpha,h}$ is a relatively compact perturbation of $\cL_{0,h} =
\LL_h -\frac32$, both operators have the same domain, and it follows
that $\omega_h$ belongs to the domain of $\LL_h$. In particular, we
have $\nabla_h \omega_h \in L^2(\infty)^4$ and $|x_h| \omega_h \in
L^2(\infty)^2$, see e.g. \cite[Section~2]{GW3}.

In the rest of the proof, we denote by $\langle \cdot,\cdot\rangle$
the inner product in the complexified space $L^2(\infty)^2$, 
namely
\[
  \langle \omega_h^1,\omega_h^2\rangle \,=\, \int_{\R^2} 
  p(x_h)\omega_h^1(x_h)\cdot\overline{\omega_h^2(x_h)}\dd x_h~,
\]
where $p(x_h) = \rho_\infty(|x_h|^2) = e^{|x_h|^2/4}$. We also 
denote $\|\omega_h\|^2 = \langle \omega_h,\omega_h\rangle$. 
We recall that $\LL_h$ is a selfadjoint operator in $L^2(\infty)^2$
which satisfies $-\LL_h \ge 0$ on $L^2(\infty)^2$ and $-\LL_h \ge 1/2$ 
on $L^2_0(\infty)^2$. For later use, we observe that the (unbounded)
operator $\omega_h \mapsto (U^G_h,\nabla_h)\omega_h$ is skew-symmetric
in $L^2(\infty)^2$, because the vector field $p(x_h)U^G(x_h)$ is 
divergence-free. 

We now take the inner product of \eqref{eq.eigenvalue2} with $\omega_h$, 
and evaluate the real part of the result. Using the skew-symmetry of 
the operator $(U^G_h,\nabla_h)$, we easily obtain
\begin{align}\label{eq.spec1}
  &\Re\lambda\,\|\omega_h\|^2 \,=\, \langle \LL_h \omega_h,\omega_h\rangle
  -\frac32\|\omega_h\|^2 + \alpha \Re\langle (\omega_h,\nabla_h)U^G_h,
  \omega_h\rangle \\ \nonumber
  \,&=\, \langle \LL_h \omega_h,\omega_h\rangle -\frac32\|\omega_h\|^2 
  + 2\alpha \Re\int_{\R^2}p(x_h)(x_h\cdot\omega_h)(x_h^\bot\cdot
  \overline{\omega_h})(u^g)'(|x_h|^2)\dd x_h~,
\end{align}
where $u^g(r)$ is defined in \eqref{def.u^g}. On the other hand, 
it follows from \eqref{eq.eigenvalue2} that the scalar function
$x_h\cdot\omega_h \in L^2(\infty)$ satisfies
\[
  \lambda\,x_h\cdot\omega_h \,=\, \LL_h(x_h\cdot\omega_h) 
  -2 x_h\cdot\omega_h -\alpha (U^G_h,\nabla_h)(x_h\cdot\omega_h) 
  -2\nabla_h \cdot \omega_h~.
\]
Thus, proceeding as above and using the same notation $\langle \cdot,
\cdot\rangle$ for the inner product in $L^2(\infty)$, we find
\begin{equation}\label{eq.spec2}
  \Re\lambda\,\|x_h\cdot\omega_h\|^2 \,=\, 
  \langle\LL_h (x_h\cdot\omega_h),x_h\cdot\omega_h \rangle
  -2 \|x_h\cdot\omega_h\|^2 -2\Re \langle \nabla_h\cdot\omega_h,
 x_h\cdot\omega_h\rangle~.
\end{equation}
Finally, the two-dimensional divergence $\nabla_h \cdot \omega_h 
\in L^2_0(\infty)$ satisfies
\begin{equation}\label{eq.2Ddiv}
  \lambda\,\nabla_h\cdot\omega_h \,=\, \LL_h(\nabla_h\cdot\omega_h) 
  - \nabla_h\cdot\omega_h -\alpha (U^G_h,\nabla_h)(\nabla_h\cdot
  \omega_h)~,
\end{equation}
hence
\begin{equation}\label{eq.spec3}
  \Re\lambda\,\|\nabla_h\cdot\omega_h\|^2 \,=\, 
  \langle\LL_h (\nabla_h\cdot\omega_h),\nabla_h\cdot\omega_h \rangle
  - \|\nabla_h\cdot\omega_h\|^2~.
\end{equation}

Since $\nabla_h\cdot\omega_h \in L^2_0(\infty)$, it follows from 
\eqref{eq.spec3} that $\Re\lambda\,\|\nabla_h\cdot\omega_h\|^2 \le 
-\frac32 \|\nabla_h\cdot\omega_h\|^2$. Thus we must have $\Re\lambda
\le -\frac32$, unless $\nabla_h\cdot\omega_h \equiv 0$. In the 
latter case, we deduce from \eqref{eq.spec2} that $\Re\lambda\,
\|x_h\cdot\omega_h\|^2 \le -2\|x_h\cdot\omega_h\|^2$, hence 
$\Re\lambda \le -2$ unless $x_h\cdot\omega_h \equiv 0$. But 
if this last condition is met, it follows from \eqref{eq.spec1}
that $\Re\lambda\,\|\omega_h\|^2 \le -\frac32 \|\omega_h\|^2$, 
hence $\Re\lambda \le -\frac32$ because $\omega_h$ is not 
identically zero. Summarizing, we conclude that $\Re \lambda \le 
-\frac32$ in all cases. \QED

\medskip\noindent{\bf Remark.} Actually the conclusions of
Proposition~\ref{prop.discrete.spectrum.2} can be slightly
strengthened. First, in the invariant subspace where $\nabla_h \cdot
\omega_h = 0$, one can show that all eigenvalues of $\cL_{\alpha,h}$
satisfy $\Re\lambda \le -2$. This follows from the proof above if we
use in addition the fact that $\omega_h \in L^2_0(\infty)^2$, due to
the divergence-free condition. The result is clearly sharp, because if
$g(x_h)$ is defined by \eqref{def.g} it is easy to verify that the function
$\omega_h = x_h^\bot g(x_h)$ satisfies $\cL_{\alpha,h}\omega_h =
-2\omega_h$ for any $\alpha \in \R$. On the other hand, if $\omega_h$
is a solution of \eqref{eq.eigenvalue2} such that $\nabla_h \cdot
\omega_h \neq 0$, we have $\Re\lambda < -\frac32$ if $\alpha \neq
0$. This follows from \eqref{eq.spec3}, because we know from
\cite[Appendix~A]{GW0} that
\[
  \langle\LL_h (\nabla_h\cdot\omega_h),\nabla_h\cdot\omega_h
  \rangle \,<\, -\frac12\|\nabla_h\cdot\omega_h\|^2~,
\]
unless $\nabla_h \cdot\omega_h = (a_1x_1 + a_2x_2)g(x_h)$ for some
$a_1, a_2 \in \C$. But this ansatz is not compatible with
\eqref{eq.2Ddiv} if $\alpha \neq 0$. In fact, using the techniques
developped in \cite{M3} or \cite{GGN}, it is possible to show that,
given any $M > 0$, the eigenvalue equation \eqref{eq.2Ddiv} restricted
to the orthogonal complement of the space of all radially symmetric
functions in $L^2(\infty)$ has no nontrivial solution such that $\Re
\lambda \ge -M$, if $|\alpha|$ is sufficiently large depending on $M$.

\medskip It is now easy to conclude the proof of
Proposition~\ref{prop.vectorial.2D}.  As was already mentioned, we
only need to prove that estimate \eqref{est.prop.vectorial.2D.1} holds
for any $\mu < 3/2$.  If $\rho_\alpha(m) > 0$ denotes the spectral radius
of the operator $e^{\cL_{\alpha,h}}$ in $L^2(m)^2$, this is equivalent
to showing that $\log \rho_\alpha(m) \le -3/2$, see
\cite[Proposition~IV.2.2]{EN}.  But that inequality follows
immediately from Propositions~\ref{prop.essential.spectrum},
\ref{prop.discrete.spectrum}, and \ref{prop.discrete.spectrum.2},
since $m > 1$. The proof of Proposition~\ref{prop.vectorial.2D} is now
complete. \QED


\section{Linear stability}\label{sec.linear}

Equipped with the results of the previous section, we now study the 
linearized equation \eqref{eq.linear} in its full generality. 
Using Proposition~\ref{prop.semigroup.L} and a perturbation argument, 
it is not difficult to verify that the linear operator $L - 
\alpha\Lambda$ generates a locally bounded semigroup in the space
$\X(m)$ for any $\alpha \in \R$ and any $m \in (1,\infty]$, see 
Proposition~\ref{prop.linear.local} below. The goal of this 
section is to show that the semigroup $e^{t(L-\alpha\Lambda)}$ extends
to a bounded operator from $\X^p(m)$ to $\X(m)$ for any $t > 0$ and
any $p \in [1,2]$, and satisfies the following uniform estimates: 

\begin{prop}\label{prop.linear} Fix $m \in (1,\infty]$, $p\in [1,2]$, 
$\alpha\in \R$, and take $\mu \in (1,\frac32)$, $\eta \in (0,\frac12]$
such that $2\mu < m+1$ and $2\eta < m-1$. For any $\beta = (\beta_1,
\beta_2,\beta_3) \in \N^3$, there exists $C > 0$ such that
\begin{align}
  \|\partial_x^\beta(e^{t(L-\alpha\Lambda)}\omega_0)_h\|_{X(m)^2} 
  \,&\le\, \frac{C\,e^{-(\mu+\beta_3)t}}{a(t)^{\frac{1}{p}-\frac{1}{2}
  +\frac{|\beta|}{2}}}\|\omega_0\|_{\X^p(m)}~,
  \label{est.prop.linear.1}\\
  \|\partial_x^\beta(e^{t(L-\alpha\Lambda)}\omega_0)_3\|_{X(m)} 
  \,&\le\, \frac{C\, e^{-(\eta+\beta_3)t}}{a(t)^{\frac{1}{p}-\frac{1}{2}
  +\frac{|\beta|}{2}}}\|\omega_0\|_{\X^p(m)}~, 
  \label{est.prop.linear.2}
\end{align}
for any $\omega_0\in \X^p(m)$ and all $t > 0$. Moreover, $\nabla\cdot 
\omega_0 = 0$, then $\nabla\cdot e^{t(L-\alpha\Lambda)}\omega_0 = 0$ 
for all $t>0$. 
\end{prop}

The proof of this important result is divided into several steps. 

\subsection{Global existence and short time estimates}\label{subsec.est.local}

We first prove that the linearized equation \eqref{eq.linear}
has a unique global solution in $\X(m)$. 

\begin{prop}\label{prop.linear.local} Fix $m\in (1,\infty]$, $p\in [1,2]$, 
and $\alpha\in \R$. Then, for any $\omega_0\in \X^p (m)$, 
Eq.~\eqref{eq.linear} has a unique solution $\omega\in L^\infty_{loc}
(\R_+;\X(m))\cap C([0,\infty);\X_{loc}^p(m))$ with initial data 
$\omega_0$. Moreover, for any $\beta \in \N^3$, there exist positive 
constants $C_1$, $C_2$ (independent of $\alpha$) such that
\begin{equation}
  \|\partial_x^\beta\omega(t)\|_{\X(m)}\le \frac{C_1}{a(t)^{\frac{1}{p}
  -\frac{1}{2}+\frac{|\beta|}{2}}}\|\omega_0\|_{\X^p(m)}~, \qquad
  \hbox{for} \quad 0 < t\le \frac{C_2}{|\alpha|^2+1}~,
  \label{est.prop.linear.local}
\end{equation}
where $a(t) = 1-e^{-t}$. Finally, if $\nabla\cdot \omega_0=0$, then 
$\nabla\cdot\omega(t)=0$ for all $t>0$.
\end{prop}

\noindent{\bf Proof.} We proceed as in the proof of 
Lemma~\ref{prop.vectorial.2D.local}. Let $e^{tL}$ be the semigroup 
generated by $L$, which is given by the explicit expression 
\eqref{eq.semigroup.L_3}. The integral equation corresponding to
\eqref{eq.linear} is
\begin{equation}
  \omega (t) \,=\, e^{tL}\omega_0 - \alpha \int_0^t e^{(t-s)L}\Lambda
  \omega(s) \dd s  \,=:\, (F\omega)(t)~, \quad t > 0~.
  \label{eq.integral.prop.linear.local}
\end{equation}
Given $k \in \N\setminus\{0\}$ and a sufficiently small $T \in (0,1]$, 
we shall solve \eqref{eq.integral.prop.linear.local} in the Banach 
space
\[
  \U_{k,T} \,=\, \Bigl\{\omega \in L^\infty_{loc}((0,T);\X(m))\cap 
  C([0,T];\X^p_{loc}(m))~\Big|~ \|\omega\|_{k,T} < \infty\Bigr\}~,
\]
equipped with the norm
\[  
  \|\omega\|_{k,T} \,=\, \sum_{|\beta|\le k} \Bigl(
  \sup_{0<t<T}a(t)^{\frac{1}{p}-\frac{1}{2}+\frac{|\beta|}{2}}
  \|\partial_x^\beta \omega(t)\|_{\X(m)} + \sup_{0<t<T}
  a(t)^{\frac{|\beta|}{2}}\|\partial_x^\beta \omega(t)\|_{\X^p(m)} 
  \Bigr)~,
\]
where $a(t) = 1 - e^{-t}$. If $\omega_0 \in \X^p(m)$, we know from 
Proposition~\ref{prop.semigroup.L} that the map $t \mapsto 
e^{tL}\omega_0$ belongs to $\U_{k,T}$ for any $T > 0$, and that
$\|e^{tL}\omega_0\|_{k,T} \le C_1\|\omega_0\|_{\X^p(m)}$ for some
$C_1 > 0$ depending only on $k$, $m$, $p$. 

Given $\omega \in \U_{k,T}$, we now estimate the integral term 
in \eqref{eq.integral.prop.linear.local}. Using 
Proposition~\ref{prop.semigroup.L} and Lemma~\ref{lem.after.Riesz}, 
we find
\begin{align}
  \|\partial_x^\beta e^{(t-s)L}\Lambda \omega(s)\|_{\X(m)} \,&\le\, 
  \frac{C\|\Lambda\omega(s)\|_{\X^p(m)}}{a(t-s)^{\frac{1}{p}-\frac{1}{2}
  +\frac{|\beta|}{2}}} \,\le\, \frac{C\sum_{|\tilde \beta|\le 1}
  \|\partial_x^{\tilde\beta}\omega(s)\|_{\X^p (m)}}{a(t-s)^{\frac{1}{p}
  -\frac{1}{2}+\frac{|\beta|}{2}}}\nonumber\\
  \,&\le\, \frac{C\|\omega\|_{k,T}}{a(t-s)^{\frac{1}{p}-\frac{1}{2}
  +\frac{|\beta|}{2}}\,a(s)^{\frac{1}{2}}}~, \qquad 0 < s < t~.
  \label{est.proof.prop.linear.local.1}
\end{align}
Similarly we have $\|\partial_x^\beta e^{(t-s)L}\Lambda \omega(s)
\|_{\X^p (m)} \le C a(t-s)^{-\frac{|\beta|}{2}}\,a(s)^{-\frac12}
\|\omega\|_{k,T}$ for $0 < s < t$. In the particular case 
where $\beta = 0$, it follows that
\begin{align}
  &\Bigl\|\int_0^t e^{(t-s)L}\Lambda \omega(s) \dd s\Bigr\|_{\X(m)}
  \,\le\, C a(t)^{1-\frac{1}{p}}\|\omega\|_{k,T}~,\label{eqint1}\\
  &\Bigl\|\int_0^t e^{(t-s)L}\Lambda \omega(s)\dd s\Bigr\|_{\X^p(m)}
  \,\le\, C a(t)^{\frac{1}{2}}\|\omega\|_{k,T}~, \qquad
  0 < t \le T~.\label{eqint2} 
\end{align}

Assume now that $1 \le |\beta| \le k$. If $\beta' \le \beta$ and 
$|\beta'| = |\beta| - 1$, we have from Lemma~\ref{lem.after.Riesz}
\[
  \|\partial_x^{\beta'} \Lambda \omega(s)\|_{\X(m)} \,\le\, 
  C \sum_{|\tilde \beta| = |\beta|}\|\partial_x^{\tilde \beta}
  \omega(s)\|_{\X(m)} \,\le\, \frac{C \|\omega\|_{k,T}}{
  a(s)^{\frac{1}{p}-\frac{1}{2} +\frac{|\beta|}{2}}}~,
  \qquad 0 < s \le T~.
\]
Thus, writing $\partial_x^\beta e^{(t-s)L} = \partial_x^{\beta-\beta'} 
\partial_x^{\beta'} e^{(t-s)L} = \partial_x^{\beta-\beta'}
\,e^{(\frac{\beta_1'+\beta_2'}{2}-\beta_3')t}e^{(t-s)L}\partial_x^{\beta'}$, 
and using Proposition~\ref{prop.semigroup.L} again, we obtain
\begin{equation}
  \|\partial_x^\beta e^{(t-s)L}\Lambda \omega(s)\|_{\X(m)} \,\le\, 
  C \|\partial_x^{\beta-\beta'}e^{(t-s)L}\partial_x^{\beta'}\Lambda 
  \omega(s)\|_{\X(m)} \,\le\, \frac{C\|\omega\|_{k,T}}{a(t-s)^{\frac{1}{2}}
  \,a(s)^{\frac{1}{p}-\frac{1}{2} +\frac{|\beta|}{2}}}~,
  \label{est.proof.prop.linear.local.3}
\end{equation}
for $0 < s < t$. Similarly, we have $\|\partial_x^\beta e^{(t-s)L}
\Lambda \omega(s)\|_{\X^p(m)} \le a(t-s)^{-\frac{1}{2}}
\,a(s)^{-\frac{|\beta|}{2}}\|\omega\|_{k,T}$. Combining 
\eqref{est.proof.prop.linear.local.1} and 
\eqref{est.proof.prop.linear.local.3}, we obtain the following
estimate
\begin{align}
  \Bigl\|\partial_x^\beta &\int_0^t e^{(t-s)L}\Lambda \omega(s)
  \dd s\Bigr\|_{\X(m)}\nonumber \\
  \,&\le\, C\Bigl(\int_0^\frac{t}{2} a(t-s)^{-\frac{1}{p}+\frac{1}{2}
  -\frac{|\beta|}{2}}\,a(s)^{-\frac{1}{2}}\dd s + 
  \int_\frac{t}{2}^t a(t-s)^{-\frac{1}{2}} \,a(s)^{-\frac{1}{p}
  +\frac{1}{2}-\frac{|\beta|}{2}}\dd s\Bigr)\|\omega\|_{k,T}
  \label{eqint3}\\
  \,&\le\, C a(t)^{1-\frac{1}{p}-\frac{|\beta|}{2}} \|\omega\|_{k,T}~,
  \qquad 0 < t \le T~, \nonumber
\end{align}
which generalizes \eqref{eqint1}. Similarly, the generalization
of \eqref{eqint2} is
\begin{equation}
  \Bigl\|\partial_x^\beta \int_0^t e^{(t-s)L}\Lambda \omega(s)
  \dd s\Bigr\|_{\X^p(m)} \,\le\, C a(t)^{\frac{1}{2}-\frac{|\beta|}{2}} 
  \|\omega\|_{k,T}~, \qquad 0 < t \le T~. \label{eqint4} 
\end{equation}

Summarizing, we have shown that the linear map $F$ defined by 
\eqref{eq.integral.prop.linear.local} satisfies the estimate
\begin{equation*}
  \|F\omega\|_{k,T} \,\le\, C_1\|\omega_0\|_{\X^p(m)} + 
  \tilde C|\alpha| T^\frac{1}{2}\|\omega\|_{k,T}~, \qquad 
  \hbox{if} \quad 0 < T \le 1~,
\end{equation*}
where $\tilde C > 0$ depends only on $k$, $m$ and $p$. Arguing
as in \cite[Corollary A.7 and Remark A.8]{GW2}, it is also straightforward 
to verify that $F\omega\in C([0,T];\X_{loc}^p(m))$ if $\omega \in 
\U_{k,T}$. If we now assume that $T \le C_2(1+|\alpha|^2)^{-1}$, where
$C_2 = 1/(4\tilde C^2)$, we see that $F$ is a strict contraction
in $\U_{k,T}$. As a consequence, the integral equation 
\eqref{eq.integral.prop.linear.local} has a unique fixed point 
$\omega \in \U_{k,T}$, which satisfies $\|\omega\|_{k,T} \le
2C_1\|\omega_0\|_{\X^p(m)}$. This proves that equation \eqref{eq.linear}
is locally well-posed in $\X^p(m)$, and since the local existence
time $T$ is independent of the initial data, the solutions can 
be extended globally in time. Finally, since both operators 
$L$ and $\Lambda$ preserve the divergence-free condition, it 
is easy to check that, if $\nabla\cdot \omega_0=0$, then 
the solution $\omega$ of \eqref{eq.linear} satisfies $\nabla \cdot 
\omega(t)=0$ for all $t > 0$. This completes the proof. \QED

\subsection{Decay estimates for the vertical derivatives}\label{subsec.est.x_3}

Proposition~\ref{prop.linear.local} shows that the linearized 
equation \eqref{eq.linear} is globally well-posed in the space
$\X(m)$ for $m > 1$, but does not provide accurate estimates 
on the solution $\omega(t) = e^{t(L-\alpha\Lambda)}\omega_0$ 
for large times. In this section, we focus on the derivatives of
$\omega(t)$ with respect to the vertical variable $x_3$. 
Using identity \eqref{eq.semigroup}, we shall show that 
$\partial_{x_3}^k \omega(t)$ decays exponentially as $t \to \infty$, 
provided $k \in \N$ is large enough depending on $|\alpha|$. 
Albeit elementary, this observation plays a crucial role in the 
proof of Proposition~\ref{prop.linear}, because it will allow
us to simplify the study of the semigroup $e^{t(L-\alpha\Lambda)}$
by disregarding most of the terms involving a vertical derivative. 

\begin{prop}\label{prop.derivative.x_3} Fix $m\in (1,\infty]$.
There exist positive constants $C_3, C_4$ such that, for all 
$\alpha \in \R$, all $k \in \N$, and all $\omega_0 \in \X(m)$ 
with $\partial_{x_3}^k \omega_0 \in \X(m)$, the following estimate 
holds:
\begin{equation}
  \|\partial_{x_3}^k e^{t(L-\alpha\Lambda)}\omega_0\|_{\X(m)}\,\le\, 
  C_3 \,e^{(C_4(|\alpha|^2+1) - k)t} \|\partial_{x_3}^k \omega_0\|_{\X(m)}~,
  \qquad t \ge 0~. \label{est.deriv.x_3}
\end{equation}

\end{prop}

\noindent{\bf Proof.} In view of \eqref{eq.semigroup}, it is 
sufficient to prove \eqref{est.deriv.x_3} for $k = 0$. If 
$\omega_0 \in \X(m)$, we know from Proposition~\ref{prop.linear.local}
that there exist constants $C_1 \ge 1$ and $C_2 > 0$, depending 
only on $m$, such that the solution $\omega(t) = e^{t(L-\alpha\Lambda)}
\omega_0$ of \eqref{eq.linear} satisfies $\|\omega(t)\|_{\X(m)}
\le C_1 \|\omega_0\|_{\X(m)}$ for $t \in (0,t_0]$, where $t_0 = 
C_2/(|\alpha|^2+1)$. Using the semigroup property, we can iterate
this bound, and we easily obtain
\[
  \|e^{t(L-\alpha\Lambda)}\omega_0\|_{\X(m)} \,\le\, C_3 
  \,e^{C_4(|\alpha|^2+1)t} \|\omega_0\|_{\X(m)}~, \qquad t \ge 0~,
\]
where $C_3 = C_1$ and $C_4 = C_2^{-1}\log(C_1)$. This concludes 
the proof. \QED

\subsection{Decomposition of the linearized operator}\label{subsec.decomposition}

Motivated by Proposition~\ref{prop.derivative.x_3}, we now decompose
the linear operator $L - \alpha\Lambda$ as follows:
\begin{equation}
  L-\alpha\Lambda \,=\, \cL_\alpha + \LL_3 - \alpha H~,
  \label{def.decompose}
\end{equation}
where $\cL_\alpha$ is defined in \eqref{def.calLambda} and $\LL_3 = 
\partial_{x_3}^2 - x_3\partial_{x_3}$. We recall that the operator
$\cL_\alpha$ does not involve any derivative with respect to the
vertical variable $x_3$, and does not couple the horizontal and
vertical components of $\omega = (\omega_h,\omega_3)^\top$. In view of
\eqref{eq.Lambda_dec}--\eqref{def.calLambda}, the last term in
\eqref{def.decompose} has the following expression:
\[
   H \,=\, \Lambda_3 - \tilde \Lambda_3 -\Lambda_4~, 
\]
where $\Lambda_3, \Lambda_4$ are defined in \eqref{def.Lambda_i} and 
$\tilde \Lambda_3$ after \eqref{def.calLambda}. More explicitly, 
we have
\begin{equation}
 H\omega \,=\, \begin{pmatrix} H_h\omega \\ H_3\omega\end{pmatrix}
 \,=\, \begin{pmatrix} -g (K_{3D}*\partial_{x_3}\omega)_h \\
 (K_{3D}*\omega-K_{2D}\star\omega_3,\nabla )g - g(K_{3D}*\partial_{x_3}
 \omega )_3 \end{pmatrix}~, \label{def.Hexp}
\end{equation}
where $K_{3D}$, $K_{2D}$ are the Biot-Savart kernels \eqref{def.3DBS}, 
\eqref{def.2DBS}, and $g$ is defined in \eqref{def.g}. Here $\star$ 
denotes the convolution with respect to the horizontal variables, 
so that
\[
  (K_{2D}\star\omega_3)(x_h,x_3) \,=\, \int_{\R^2}K_{2D}(x_h-y_h)
  \,\omega_3(y_h,x_3)\dd y_h~.
\]
Thus, unlike $\cL_\alpha$, the operator $H$ involves vertical 
derivatives, and couples the horizontal and vertical components of 
$\omega$. As was already observed in Section~\ref{sec.vectorial}, 
we have $H\omega = 0$ whenever $\partial_{x_3}\omega = 0$, see 
Proposition~\ref{prop.H} below.

Let $R_\alpha (t)$ denote the semigroup generated by the linear 
operator $\cL_\alpha + \LL_3$. In analogy with \eqref{eq.semigroup.L_3}, 
we have the following representation:
\begin{equation}
  (R_\alpha(t)\omega)(x) \,=\, \frac{1}{\sqrt{2\pi a(2t)}}\int_{\R}
  e^{-\frac{|x_3e^{-t}-y_3|^2}{2a(2t)}} \Bigl(e^{t\cL_\alpha}\omega
  (\cdot,y_3)\Bigr)(x_h)\dd y_3~, \quad t > 0~,
  \label{eq.semigroup.R_alpha}
\end{equation}
where $a(t) = 1-e^{-t}$ and $e^{t\cL_\alpha}$ is the semigroup 
generated by $\cL_\alpha$. Since $R_\alpha (t)$ does not couple
the horizontal and vertical components of $\omega$, we can 
write
\[
  R_{\alpha}(t)\omega \,=\, \begin{pmatrix}
  R_{\alpha,h}(t)\omega_h \\ R_{\alpha,3}(t)\omega_3 
 \end{pmatrix}~,
\]
where $R_{\alpha,h}(t)$ and $R_{\alpha,3}(t)$ are the semigroups
generated by $\cL_{\alpha,h} + \LL_3$ and $\cL_{\alpha,3} + \LL_3$, 
respectively. Using the results of Section~\ref{sec.vectorial}, 
we obtain the following estimates:

\begin{prop}\label{prop.R_alpha} Fix $m \in (1,\infty]$, 
$\alpha \in \R$, $\mu \in (1,\frac32)$, and take $\eta \in (0,\frac12]$ 
such that $2\eta < m-1$. Then there exists $C_5 > 0$ such that 
\begin{align}
  \|R_{\alpha,h}(t)\omega_h\|_{X(m)^2} \,&\le\, C_5\,e^{-\mu t}
  \|\omega_h\|_{X(m)^2}~, \label{est.R1}\\
  \|R_{\alpha,3}(t)\omega_3\|_{X(m)} \,&\le\, C_5\,e^{-\eta t}
  \|\omega_3\|_{X(m)}~, \label{est.R2}
\end{align}
for all $\omega \in \X(m)$ and all $t \ge 0$.
\end{prop}

\noindent{\bf Proof.} Both estimates follow from the representation 
\eqref{eq.semigroup.R_alpha}, Proposition~\ref{prop.vectorial.2D}, 
and estimate \eqref{est.semigroup_3}. The calculations are 
straightforward, and can be omitted here. We just remark that, 
even if $\nabla\cdot\omega = 0$, the map $x_h \mapsto \omega_h(x_h,x_3)$ 
usually has a nonzero divergence for all values of $x_3 \in \R$. 
This is why Proposition~\ref{prop.vectorial.2D}, hence also
Proposition~\ref{prop.R_alpha}, was established without imposing
any divergence-free condition. \QED

\medskip

We conclude this section with a useful bound on the linear 
operator $H$. 

\begin{prop}\label{prop.H} Fix $m\in (1,\infty]$ and $\gamma \in
(0,1)$. There exists $C_6 > 0$ such that, for all $\omega \in \X(m)$
with $\partial_{x_3}\omega \in \X(m)$, one has
\begin{align}
  \|H_h\omega\|_{X(m)^2} \,&\le\, C_6\|\partial_{x_3}\omega\|_{\X(m)}~,
    \label{est.H1}\\ 
  \|H_3\omega\|_{X(m)} \,&\le\, C_6(\|\partial_{x_3}\omega\|_{\X(m)} + 
  \|\omega_h\|_{X(m)^2}^\gamma \|\partial_{x_3}\omega_h\|_{X(m)^2}^{1-\gamma})~.
\label{est.H2}
\end{align}
\end{prop}

\noindent{\bf Proof.} We use the expression \eqref{def.Hexp} of the
linear operator $H$. Since $\partial_{x_3}\omega \in \X(m)$, we know
from Proposition~\ref{prop.Riesz} that $\partial_{x_3}u \equiv
K_{3D}*\partial_{x_3}\omega \in X^4(0)$. Thus, using H\"older's
inequality, we obtain
\[
\|g\,\partial_{x_3}u\|_{\X(m)} \,\le\, \|\partial_{x_3}u\|_{X^4(0)}
\Bigl(\int_{\R^2}\rho_m(|x_h|^2)^2 g(x_h)^4\dd x_h\Bigr)^{1/4}
\,\le\, C \|\partial_{x_3}\omega\|_{\X(m)}~.
\]
In particular, we have $\|H_h\omega\|_{X(m)^2} \le C\|\partial_{x_3}
\omega\|_{\X(m)}$.

We next consider the two-dimensional vector $I = (K_{3D}*\omega -
K_{2D}\star\omega_3)_h$ and estimate the term $(I,\nabla_h)g$. Using the
definitions \eqref{def.3DBS}, \eqref{def.2DBS}, it is straightforward
to verify that $I(x) = I_1(x) + I_2(x)$, where
\begin{align*}
  I_1(x) \,&=\, \frac{1}{4\pi} \int_{\R^3} \frac{(x_h{-}y_h)^\bot}{
  |x-y|^3}\,(\omega_3(y_h,y_3)-\omega_3(y_h,x_3)) \dd y~, \\ 
  I_2(x) \,&=\, \frac{1}{4\pi} \int_{\R^3} \frac{(x_3{-}y_3)}{|x-y|^3}
  \,(\omega_h(y_h,y_3)-\omega_h(y_h,x_3))^\bot \dd y~.
\end{align*}
Since $\nabla_h g (x_h)=- g(x_h) x_h/2$ and $|x_h\cdot
(x_h-y_h)^\bot|\leq |x_h||x_h-y_h|^{1-\sigma}|y_h|^\sigma$ for any
$\sigma\in [0,1]$, we can bound
\begin{align*}
  |(I_1,\nabla_h) g(x)| \,&\le\, C g(x_h) |x_h| \int_{|x_3-y_3|\ge 1}
  \frac{|y_h|^\sigma}{|x-y|^{2+\sigma}}\,|\omega_3(y_h,y_3)-
  \omega_3(y_h,x_3)|\dd y \\
  &\quad~ + C g(x_h) |x_h| \int_{|x_3-y_3|< 1}\frac{1}{|x-y|^{2}}
  \,|\omega_3(y_h,y_3)-\omega_3(y_h,x_3)| \dd y~.
\end{align*}
We now proceed like in the proof of Lemma~\ref{lem.Riesz}. Integrating
first with respect to the horizontal variable $y_h \in \R^2$ and applying 
H\"older's inequality, we obtain
\begin{align*}
 |(I_1,\nabla_h) g(x)| \,&\le\, C g(x_h) |x_h| \int_{|x_3-y_3|\ge 1}
  \frac{1}{|x_3-y_3|^{2+\sigma}} \||\cdot|^\sigma \{\omega_3(\cdot,y_3)-
  \omega_3(\cdot,x_3)\}\|_{L^1(\R^2)} \dd y_3 \\
  &\quad~  + C g(x_h) |x_h| \int_{|x_3-y_3|< 1}\frac{1}{|x_3-y_3|} 
  \|\omega_3(\cdot,y_3)-\omega_3(\cdot,x_3) \|_{L^2(\R^2)} 
  \dd y_3~.
\end{align*}
Assuming $0 < \sigma < m-1$, we have the estimate $\||\cdot|^\sigma f
\|_{L^1(\R^2)} \le C\|f\|_{L^2(m)}$ for any $f \in L^2(m)$, hence
\[
 \||\cdot|^{\sigma} \{\omega_3 (\cdot,y_3)-\omega_3 (\cdot,x_3)\}
 \|_{L^1(\R^2) } + \|\omega_3(\cdot,y_3)-\omega_3(\cdot,x_3)\|_{L^2(\R^2)} 
 \,\le\, C |x_3-y_3| \|\partial_{x_3}\omega_3\|_{X(m)}~.
\]
We conclude that
\begin{align*}
  |(I_1,\nabla_h) g(x)| \,&\le\, C g(x_h)|x_h| \Bigl(\int_{|x_3-y_3|\ge 1}
  \frac{1}{|x_3-y_3|^{1+\sigma}} \dd y_3 + \int_{|x_3-y_3|< 1} \dd y_3 \Bigr)
  \|\partial_{x_3}\omega_3\|_{X(m)}\\
  \,&\le\, C g(x_h)|x_h|\|\partial_{x_3}\omega_3\|_{X(m)}~,
\end{align*}
which gives the bound $\|(I_1,\nabla)g\|_{X(m)}\le C\|\partial_{x_3}
\omega_3\|_{X(m)}$.  

Finally we consider the term $(I_2,\nabla_h)g$. Using again H\"older's 
inequality, we obtain
\begin{align*}
  |I_2(x)| \,&\le\, C\int_{|x_3-y_3|\ge 1}\frac{1}{|x-y|^2}|\omega_h
  (y_h,y_3)-\omega_h(y_h,x_3)|\dd y \\
  &\quad~ + C \int_{|x_3-y_3|< 1}\frac{1}{|x-y|^2}|\omega_h (y_h,y_3)-
  \omega_h (y_h,x_3)|\dd y \\
  \,&\le\, C \int_{|x_3-y_3|\ge 1}\frac{1}{|x_3-y_3|^2}\|\omega_h(\cdot,y_3)-
  \omega_h(\cdot,x_3)\|_{L^1(\R^2)} \dd y_3 \\
  &\quad~ +  \int_{|x_3-y_3|< 1}\frac{1}{|x_3-y_3|}\|\omega_h(\cdot,y_3)-
  \omega_h(\cdot,x_3) \|_{L^2(\R^2)}\dd y_3~.
\end{align*}
Since $L^2(m) \hookrightarrow L^p(\R^2)$ for $p \in [1,2]$, we have
$\|\omega_h (\cdot,y_3)-\omega_h (\cdot,x_3)\|_{L^p(\R^2)^2} \le
2\|\omega_h\|_{X(m)^2}$ and $\|\omega_h (\cdot,y_3)-\omega_h
(\cdot,x_3)\|_{L^p(\R^2)^2} \le |x_3-y_3|\|\partial_{x_3}\omega_h\|_{X(m)^2}$.
In particular, for any $\gamma \in (0,1)$,
\[
  \|\omega_h (\cdot,y_3)-\omega_h (\cdot,x_3)\|_{L^p(\R^2)^2} \,\le\, 
  2^\gamma |x_3-y_3|^{1-\gamma} \|\omega_h \|_{X(m)^2}^\gamma 
  \|\partial_{x_3}\omega_h \|_{X(m)^2}^{1-\gamma}~.
\]
Thus we obtain
\[
  \|I_2\|_{L^\infty(\R^3)^2} \,\le\, C \|\omega_h\|_{X(m)^2}^\gamma 
  \|\partial_{x_3}\omega_h\|_{X(m)^2}^{1-\gamma}~,
\]
and conclude that $\|(I_2,\nabla_h)g\|_{X(m)} \le
C\|\omega_h\|_{X(m)^2}^\gamma \|\partial_{x_3}\omega_h\|_{X(m)^2}^{1-\gamma}$.
This completes the proof of Proposition \ref{prop.H}. \QED

\subsection{Large time estimates}\label{subsec.est.large}

In this section we complete the proof of Proposition~\ref{prop.linear}. 
Fix $m \in (1,\infty]$, $\alpha \in \R$, and assume that $\omega_0 
\in \X^p(m)$ for some $p \in [1,2]$. Let $\omega(t) = e^{t(L-\alpha\Lambda)}
\omega_0$ be the solution of the linearized equation \eqref{eq.linear}
given by Proposition~\ref{prop.linear.local}. Take any $k \in \N$ 
such that $k > C_4(|\alpha|^2+1) + 1/2$, where $C_4$ is as in 
Proposition~\ref{prop.derivative.x_3}, and choose $t_0 > 0$ 
small enough so that estimate \eqref{est.prop.linear.local} holds 
for all $t \in (0,t_0]$ and all $\beta \in \N^3$ with $|\beta| \le k$. 
Our goal is to control the solution $\omega(t)$ for $t \ge t_0$ and
to establish the decay estimates \eqref{est.prop.linear.1}, 
\eqref{est.prop.linear.2}. 

To this end, we first observe that $\omega(t)$ satisfies the 
integral equation
\begin{equation}
  \omega(t) \,=\,R_\alpha (t-t_0)\omega(t_0) - \alpha\int_{t_0}^t 
  R_\alpha (t-s) H\omega(s)\dd s~, \quad t \ge t_0~,
  \label{eq.linear.integral}
\end{equation}
where $R_\alpha(t)$ is the semigroup defined by 
\eqref{eq.semigroup.R_alpha}. Fix $\bar\eta \in (0,1/2)$ such that
$2\bar\eta < m-1$. By Proposition~\ref{prop.R_alpha}, we have
\begin{equation}\label{est.om1}
  \|\omega(t)\|_{\X(m)} \,\le\, C_5\,e^{-\bar\eta(t-t_0)}\|\omega(t_0)\|_{\X(m)}
  + C_5|\alpha|\int_{t_0}^t e^{-\bar\eta(t-s)}\|H\omega(s)\|_{\X(m)}\dd s~.
\end{equation}
To estimate the term $\|H\omega(s)\|_{\X(m)}$, we first apply 
Proposition~\ref{prop.H} with $\gamma = 1/2$, and then the 
classical interpolation inequality
\[
  \|\partial_{x_3}\omega\|_{\X(m)} \,\le\, C \|\omega\|_{\X(m)}^{1-1/k}
  \,\|\partial_{x_3}^k\omega\|_{\X(m)}^{1/k}~.
\]
Using in addition Young's inequality, we conclude that, given any 
$\epsilon > 0$, there exists $C_\epsilon > 0$ such that
\begin{equation}\label{est.om2}
  C_5 |\alpha|\,\|H\omega(s)\|_{\X(m)} \,\le\, \epsilon \|\omega(s)\|_{\X(m)}
  + C_\epsilon \|\partial_{x_3}^k\omega(s)\|_{\X(m)}~.
\end{equation}
On the other hand, since $k > C_4(|\alpha|^2+1) + 1/2$, it follows
from \eqref{est.deriv.x_3} that
\begin{equation}\label{est.om3}
  \|\partial_{x_3}^k\omega(s)\|_{\X(m)} \,\le\, C_3\,e^{-(s-t_0)/2}
  \|\partial_{x_3}^k\omega(t_0)\|_{\X(m)}~, \quad s \ge t_0~.
\end{equation}
Replacing \eqref{est.om2} and \eqref{est.om3} into \eqref{est.om1}, 
we easily obtain
\[
  \|\omega(t)\|_{\X(m)} \,\le\, \Bigl(C_5\|\omega(t_0)\|_{\X(m)} + 
  C_\epsilon'\|\partial_{x_3}^k\omega(t_0)\|_{\X(m)}\Bigr)
  \,e^{-\bar\eta(t-t_0)} + \epsilon \int_{t_0}^t e^{-\bar\eta(t-s)}
  \|\omega(s)\|_{\X(m)}\dd s~,
\]
for some $C_\epsilon' > 0$. Applying now Gronwall's lemma, and using
\eqref{est.prop.linear.local} to bound $\|\omega(t_0)\|_{\X(m)}$ 
and $\|\partial_{x_3}^k\omega(t_0)\|_{\X(m)}$ in terms of $\omega_0$, 
we see that $\|\omega(t)\|_{\X(m)} \le C\,e^{-\eta t}\|\omega_0\|_{\X^p(m)}$
for $t \ge t_0$, where $\eta = \bar \eta - \epsilon$. Finally, 
using \eqref{est.prop.linear.local} again to control the solution
for $t < t_0$, we conclude that there exists $C_7 > 0$ such that
\begin{equation}\label{est.om4}
  \|\omega(t)\|_{\X(m)} \,\equiv\, \|e^{t(L-\alpha\Lambda)}\omega_0\|_{\X(m)}
  \,\le\, \frac{C_7\,e^{-\eta t}}{a(t)^{\frac1p-\frac12}}
  \,\|\omega_0\|_{\X^p(m)}~,
\end{equation}
for all $t > 0$. Since $\epsilon > 0$ was arbitrary, estimate 
\eqref{est.om4} holds for any  $\eta \in (0,1/2)$ such that 
$2\eta < m-1$. 

To conclude the proof, it remains to find the optimal decay rates for
$\|\omega_h(t)\|_{\X(m)}$, $\|\omega_3(t)\|_{\X(m)}$ (including the
value $\eta = 1/2$ if $m > 2$), and to establish
\eqref{est.prop.linear.1}, \eqref{est.prop.linear.2} for $\beta \neq
0$ too. First, combining \eqref{eq.semigroup}, \eqref{est.om4} and
using \eqref{est.prop.linear.local} again for short times, we easily
obtain
\begin{equation}\label{est.om5}
  \|\partial_{x_3}\omega(t)\|_{\X(m)} \,\equiv\, \|\partial_{x_3} 
  e^{t(L-\alpha\Lambda)}\omega_0\|_{\X(m)} \,\le\, \frac{C\,e^{-(\eta+1)t}}
  {a(t)^{\frac1p}}\,\|\omega_0\|_{\X^p(m)}~,
\end{equation}
for all $t > 0$. Moreover, if $m > 2$, we know from
Proposition~\ref{prop.R_alpha} that \eqref{est.om1} holds with $\bar
\eta = 1/2$. Thus, applying Proposition~\ref{prop.H} to estimate
$\|H\omega(s)\|_{\X(m)}$ and using \eqref{est.om4}, \eqref{est.om5},
we find that $\|\omega(t)\|_{\X(m)}$ decays like $e^{-t/2}$ as $t \to
\infty$, hence \eqref{est.om4} holds with $\eta = 1/2$ if $m >
2$. 

Next, to obtain a faster decay estimate for the horizontal component 
$\omega_h$, we use \eqref{est.R1} and \eqref{est.H1}. Instead of 
\eqref{est.om1}, we find
\begin{equation}\label{est.om6}
  \|\omega_h(t)\|_{X(m)^2} \,\le\, C\,e^{-\mu(t-t_0)}
  \|(\omega(t_0))_h\|_{X(m)^2} + C|\alpha|\int_{t_0}^t e^{-\mu(t-s)}
  \|\partial_{x_3}\omega(s)\|_{\X(m)}\dd s~,
\end{equation}
for any $\mu \in (1,\frac32)$. Since $\|\partial_{x_3}\omega(t)\|_{\X(m)}
\le Ce^{-(\eta+1)t} \|\partial_{x_3}\omega_0\|_{\X(m)}$ by 
\eqref{eq.semigroup}, \eqref{est.om4}, we conclude that 
$\|\omega_h(t)\|_{X(m)^2}$ decays like $e^{-\mu t}$ as $t \to \infty$, 
provided $\mu < 1+\eta$. In other words, if $\mu \in (1,\frac32)$ 
satisfies $2\mu < m+1$, we have
\begin{equation}\label{est.om7}
  \|\omega_h(t)\|_{X(m)^2} \,\equiv\, \|(e^{t(L-\alpha\Lambda)}\omega_0)_h
  \|_{X(m)^2} \,\le\, C\,e^{-\mu t}(\|(\omega_0)_h\|_{X(m)^2} 
  + \|\partial_{x_3}\omega_0\|_{\X(m)})~,
\end{equation}
for all $t > 0$. Using the arguments leading to \eqref{est.om6} 
and proceeding as in Proposition~\ref{prop.linear.local}, we can 
also derive the following short time estimate, which complements
\eqref{est.prop.linear.local}:
\begin{equation}
  \|\partial_x^\beta\omega_h(t)\|_{X(m)^2} \,\le\, \frac{C_1}{a(t)^{\frac{1}{p}
  -\frac{1}{2}+\frac{|\beta|}{2}}}\Bigl(\|(\omega_0)_h\|_{X^p(m)^2} 
  + \|\partial_{x_3}\omega_0\|_{\X^p(m)}\Bigr)~, \quad
  0 < t\le \frac{C_2}{|\alpha|^2{+}1}~. 
  \label{est.om8}
\end{equation}

Finally, to obtain decay estimates for the derivative $\partial_x^
\beta \omega(t)$, where $\beta \in \N^3$, we can restrict ourselves to 
$t \ge 2t_1$, where $t_1 > 0$ is small enough so that the short time 
estimates \eqref{est.prop.linear.local}, \eqref{est.om8} hold for
$0 < t \le 2t_1$. In view of \eqref{eq.semigroup}, we have the identity
\[
  \partial_x^\beta e^{t(L-\alpha\Lambda)}\omega_0 \,=\, 
  e^{-\beta_3(t-t_1)}~\partial_{x_h}^{\beta_h} e^{t_1(L-\alpha\Lambda)}~ 
  e^{(t-2t_1)(L-\alpha\Lambda)}~\partial_{x_3}^{\beta_3} e^{t_1(L-\alpha\Lambda)}
  \omega_0~.
\]
Using the short time estimates \eqref{est.prop.linear.local}, 
\eqref{est.om8} with $p = 2$ to bound the first operator 
$\partial_{x_h}^{\beta_h} e^{t_1(L-\alpha\Lambda)}$, then the long-time 
estimates \eqref{est.om4}, \eqref{est.om5} or \eqref{est.om7} to treat 
the middle term $e^{(t-2t_1)(L-\alpha\Lambda)}$, and finally 
\eqref{est.prop.linear.local} again to bound the last term 
$\partial_{x_3}^{\beta_3} e^{t_1(L-\alpha\Lambda)}\omega_0$, we easily 
obtain \eqref{est.prop.linear.1} and \eqref{est.prop.linear.2}, 
together with the following estimate
\begin{equation}
  \|\partial_x^\beta(e^{t(L-\alpha\Lambda)}\omega_0)_h\|_{X(m)^2} 
  \,\le\, \frac{C\,e^{-(\mu+\beta_3)t}}{a(t)^{\frac{1}{p}-\frac{1}{2}
  +\frac{|\beta|}{2}}}\Bigl(\|(\omega_0)_h\|_{X^p(m)^2} 
  + \|\partial_{x_3}\omega_0\|_{\X^p(m)}\Bigr)~, \quad t > 0, 
  \label{est.prop.linear.3}
\end{equation}
which will also be used in the next section. This concludes the proof of 
Proposition~\ref{prop.linear}. \QED


\section{Nonlinear stability}\label{sec.nonlinear}

In this section we consider the nonlinear stability of the Burgers
vortex and prove Theorems~\ref{thm.main1} and \ref{thm.main2}.  
Our starting point is the perturbation equation \eqref{eq.omega}, which 
is equivalent to the integral equation
\begin{equation}
  \omega(t) \,=\, e^{t(L-\alpha\Lambda)}\omega_0 + \sum_{j=1}^2 
  \int_0^t e^{(t-s)(L-\alpha\Lambda)} N_j(\omega (s),\omega (s))\dd s~, 
  \qquad t \ge 0~,
  \label{eq.omega.integral}
\end{equation}
where $N_1(v,w) = (K_{3D}*v,\nabla)w$, $N_2(v,w) = (v,\nabla)K_{3D}*w$, 
and $K_{3D}$ is the Biot-Savart kernel \eqref{def.3DBS}. We first 
establish the following result, which already implies 
Theorem~\ref{thm.main1}. 

\begin{prop}\label{prop.nonlinear.stability} Fix $m\in (1,\infty]$, 
$\alpha\in \R$, and take $\eta \in (0,\frac12]$ such that $2\eta < m-1$. 
Then there exist $\delta = \delta(\alpha,m,\eta) > 0$ and $C = C(\alpha,m,
\eta) > 0$ such that, for any $\omega_0\in \X(m)$ with $\nabla\cdot 
\omega_0=0$ and $\|\omega_0\|_{\X(m)} \le \delta$, 
Eq.~\eqref{eq.omega.integral} has a unique solution 
$\omega\in L^\infty(\R_+;\X(m))\cap C([0,\infty);\X_{loc}(m))$,  
which satisfies
\begin{equation}\label{eq.decay2}
   \|\partial_x^\beta \omega(t)\|_{\X(m)} \,\le\, 
   \frac{C\|\omega_0\|_{\X(m)}}{a(t)^\frac{|\beta|}{2}} 
   \,e^{-\eta t}~, \qquad t>0~,
\end{equation}
for any multi-index $\beta \in \N^3$ of length $|\beta|\le 1$.
\end{prop}

\noindent{\bf Proof.} Let $\U$ be the Banach space of all $\omega
\in L^\infty(\R_+;\X(m))\cap C([0,\infty);\X_{loc}(m))$ such that 
$\nabla\cdot \omega(t)=0$ for all $t > 0$ and $\|\omega\|_\U < \infty$, 
where
\[
  \|\omega\|_\U \,=\, \sum_{|\beta|\le 1}\sup_{t>0} a(t)^\frac{|\beta|}{2}
  e^{\eta t}\|\partial_x^\beta \omega(t)\|_{\X(m)}~.
\]
Given $\omega_0 \in \X(m)$ such that $\nabla\cdot\omega_0=0$, we denote 
by $\Phi : \U \to \U$ the nonlinear map defined by
\begin{equation}
  \Phi(\omega)(t) \,=\, e^{t(L-\alpha\Lambda)}\omega_0 + 
  \sum_{j=1}^2\Phi_j(\omega,\omega)(t)~, \quad t > 0~,
  \label{def.Phi}
\end{equation}
where $\Phi_1$, $\Phi_2$ are the following bilinear operators:
\begin{equation}
  \Phi_j(\omega,\tilde\omega)(t) \,=\, \int_0^t e^{(t-s)(L-\alpha\Lambda)}
  N_j(\omega(s),\tilde\omega(s))\dd s~, \quad j = 1,2~.
  \label{def.Phi_j}
\end{equation}
If $\|\omega_0\|_{\X(m)}$ is sufficiently small, we shall show that the 
map $\Phi$ is a strict contraction in the ball $B_K = \{\omega \in \U\,|\, 
\|\omega\|_\U \le K\}$ for some suitable $K > 0$. It will follow that 
$\Phi$ has a unique fixed point $\omega$ in $B_K$, which by construction 
is the desired solution of \eqref{eq.omega.integral}. 

Since $\omega_0 \in \X(m)$ and $\nabla\cdot\omega_0=0$, 
Proposition~\ref{prop.linear} shows that the map $t \mapsto 
e^{t(L-\alpha\Lambda)}\omega_0$ belongs to $\U$, and satisfies the 
estimate 
\begin{equation*}
   \|e^{t(L-\alpha\Lambda)}\omega_0\|_\U \,\le\, C_1\|\omega_0\|_{\X(m)}~,
\end{equation*}
for some $C_1 > 0$ (depending on $m$, $\alpha$, $\eta$). On the other
hand, if $v,w \in \X(m)$, Corollary~\ref{cor.prop.Riesz}
implies that $N_1(v,w)$ and $N_2(v,w)$ belong to $X^p(m)^3$ for any 
$p \in (1,2)$, and satisfy the bound
\[
  \|N_1(v,w)\|_{X^p(m)^3} + \|N_2(v,w)\|_{X^p(m)^3} \,\le\, 
  C\|v\|_{\X(m)}\|\nabla w\|_{\X(m)}~,
\]
for some $C > 0$ (depending on $m$ and $p$). If in addition $\nabla\cdot 
v = 0$, then denoting $u = K_{3D}*v$ we find
\begin{equation}
  \int_{\R^2} (N_1(v,v) + N_2(v,v))_3\dd x_h \,=\, \int_{\R^2} \nabla_h 
  \cdot (v_h u_3 - u_h v_3)\dd x_h \,=\, 0~,
\end{equation}
for all $x_3 \in \R$, hence $N_1(v,v) + N_2(v,v) \in \X^p(m)$. As 
a consequence, if $\omega, \tilde \omega \in \U$, we have $N_j(\omega(t),
\tilde\omega(t)) \in X^p(m)^3$ for $j = 1,2$ and all $t > 0$, and using 
Proposition~\ref{prop.linear} again we obtain the following estimate 
for the bilinear operators $\Phi_j$: 
\begin{align*}
  \Bigl\|\sum_{j=1}^2 \partial_x^\beta &\Phi_{j}(\omega,\tilde\omega)
  (t)\Bigr\|_{\X(m)} \,\le\, \sum_{j=1}^2 \int_0^t \|\partial_x^\beta 
   e^{(t-s)(L-\alpha\Lambda)}N_j(\omega(s),\tilde\omega(s))\|_{\X(m)}\dd s\\
  \,&\le\, C \sum_{j=1}^2 \int_0^t \frac{e^{-\eta(t-s)}}{a(t{-}s)^{\frac{1}{p}
  -\frac{1}{2}+\frac{|\beta|}{2}}}\|N_j(\omega(s),\tilde\omega(s))
  \|_{X^p(m)^3}\dd s\\
  \,&\le\, C \int_0^t\frac{e^{-\eta(t-s)}}{a(t{-}s)^{\frac{1}{p}-\frac{1}{2}
  +\frac{|\beta|}{2}}} \|\omega(s)\|_{\X(m)}\|\nabla\tilde 
  \omega(s)\|_{\X(m)} \dd s\\
  \,&\le\, C \int_0^t \frac{e^{-\eta (t-s)} e^{-2\eta s}}{a(t{-}s)^{\frac{1}{p}
  -\frac{1}{2}+\frac{|\beta|}{2}}a(s)^\frac{1}{2}}\dd s \,\|\omega\|_\U
  \|\tilde \omega\|_\U \,\le\, \frac{Ce^{-\eta t}}{a(t)^{\frac{1}{p}
  +\frac{|\beta|}{2}-1}}\|\omega\|_\U\|\tilde\omega\|_\U~.
\end{align*}
Since we also know that $N_1(\omega(t),\omega(t)) + N_2(\omega(t),
\omega(t))$ belongs to $\X^p(m)$ for all $t > 0$ and is divergence-free, 
we have shown that $\Phi$ maps $\U$ into $\U$, and that there exists 
$C_2 > 0$ (depending on $|\alpha|$, $m$, and $\eta$) such that
\[
  \|\Phi(\omega)\|_\U \,\le\, C_1\|\omega_0\|_{\X(m)} + C_2\|\omega\|_\U^2~, 
  \quad 
  \|\Phi(\omega) - \Phi(\tilde\omega)\|_\U \,\le\, C_2(\|\omega\|_\U +  
  \|\tilde\omega\|_\U)\|\omega-\tilde\omega\|_\U~,
\]
for all $\omega, \tilde \omega \in \U$. We now take $K > 0$ such that
$2C_2K < 1$, and assume that $\|\omega_0\|_{\X(m)} \le K/(2C_1)$. Then
the estimates above show that $\Phi$ is a strict contraction in the
ball $B_K$, hence has a unique fixed point $\omega \in B_K$ which, of 
course, satisfies \eqref{eq.omega.integral}. Moreover $\|\omega\|_\U
\le 2C_1\|\omega_0\|_{\X(m)}$, hence \eqref{eq.decay2} holds with 
$C = 2C_1$. This concludes the proof. \QED

\medskip\noindent{\bf Remark.}
The size $\delta$ of the local basin of attraction of the Burgers
vortex $\alpha G$ in $\X(m)$ depends a priori on $\alpha$, $m$, and
$\eta$. However, as announced in Theorem~\ref{thm.main2}, the
dependence on the decay rate $\eta$ can easily be removed by the
following (standard) argument. Given $m > 1$, we first choose $\eta =
\bar\eta(m) = \min(\frac12,\frac{m-1}{4})$ and apply
Proposition~\ref{prop.nonlinear.stability} with that value of
$\eta$. We thus obtain a constant $\bar\delta > 0$ depending only on
$\alpha$ and $m$ such that, for any $\omega_0\in \X(m)$ with
$\nabla\cdot \omega_0=0$ and $\|\omega_0\|_{\X(m)} \le \bar\delta$,
Eq.~\eqref{eq.omega.integral} has a unique solution $\omega \in
L^\infty(\R_+;\X(m)) \cap C([0,\infty);\X_{loc}(m))$, which
converges exponentially to zero as $t \to \infty$. In particular,
given any $\eta \in (0,\frac12]$ such that $2\eta < m-1$, there exists
$T = T(\eta) > 0$ such that $\|\omega(t)\|_{\X(m)} \le \delta$ for all
$t \ge T$, where $\delta = \delta(\alpha,m,\eta)$ is the constant
given by Proposition~\ref{prop.nonlinear.stability}. By uniqueness 
of the solution, we conclude that $\omega$ satisfies \eqref{eq.decay2} 
for any admissible value of $\eta$.

\medskip
In view of Proposition~\ref{prop.nonlinear.stability} and the
remark that follows, the proof of Theorem~\ref{thm.main2} will be 
complete once we have established the improved decay estimate 
\eqref{est.thm.main2.1} for the horizontal component $\omega_h$.
A convenient way to do so is to repeat the proof of  
Proposition~\ref{prop.nonlinear.stability} using a different 
function space, which incorporates a faster decay rate as $t \to 
\infty$. Given $\mu \in (1,1+\eta)$, where $\eta \in (0,\frac12]$ 
is as in Proposition~\ref{prop.nonlinear.stability}, we introduce
the space $\V \subset \U$ defined by the norm
\[
  \|\omega\|_\V \,=\, \sum_{k=0,1}\sum_{|\beta|\le 1}\Bigl(
  \sup_{t>0} a(t)^\frac{k}{2} e^{(\mu +k\eta)t}\|\partial_{x_3}^k
  \partial_x^\beta \omega_h (t)\|_{X(m)^2} + \sup_{t>0}a(t)^\frac{k}{2} 
  e^{(\eta+k)t}\|\partial_{x_3}^k\partial_x^\beta \omega_3(t)\|_{X(m)}
  \Bigr)~.
\]
As in the remark above, we can assume here (without loss of generality)
that $\|\partial_x^\beta\omega_0\|_{\X(m)}$ is finite and arbitrarily 
small, for all $\beta \in \N^3$ with $|\beta| \le 1$. Using 
Proposition~\ref{prop.linear}, we thus obtain
\[
  \|e^{t(L-\alpha\Lambda)}\omega_0\|_\V \,\le\, C_3\sum_{|\beta|\le 1}
  \|\partial_x^\beta\omega_0\|_{\X(m)}~,
\]
for some $C_3 > 0$. On the other hand, if $v,w \in \X(m)$, the 
following estimates hold for any $p \in (1,2)$:
\begin{align*}
  \|N_{1,h}(v,w)\|_{X^p(m)^2} \,&\le\, C\|v\|_{\X(m)}
    \|\nabla w_h\|_{X(m)^2}~, \\
  \|N_2(v,w)\|_{X^p(m)^3} \,&\le\, C(\|v_h\|_{X(m)^2}\|\nabla_h w\|_{\X(m)} + 
    C\|v_3\|_{X(m)}\|\partial_{x_3} w\|_{\X(m)})~,\\
  \|\partial_{x_3}N_j(v,w)\|_{X^p(m)^3} \,&\le\, C(\|\partial_{x_3}v\|_{\X(m)}
    \|\nabla w\|_{\X(m)} + \|v\|_{\X(m)}\|\partial_{x_3}\nabla w\|_{\X(m)})~.
\end{align*}

We now estimate the bilinear operators $\Phi_j(\omega,\tilde\omega)$ 
for $\omega, \tilde\omega \in \V$. First, using \eqref{est.prop.linear.3}, 
we find for $t \ge 1$:
\begin{align}\nonumber
  \|\partial_x^\beta &\Phi_{1,h}(\omega,\tilde\omega)(t)\|_{X(m)^2}
  \,\le\, \int_0^t \|\partial_x^\beta \{e^{(t-s)(L-\alpha\Lambda)}
  N_1(\omega(s),\tilde\omega(s))\}_h\|_{X(m)^2}\dd s\\ \nonumber
  \,&\le\, C \int_0^t \frac{e^{-\mu(t-s)}}{a(t{-}s)^{\frac{1}{p}
  -\frac{1}{2}+\frac{|\beta|}{2}}} (\|N_{1,h}(\omega(s),\tilde
  \omega(s))\|_{X^p(m)^2} + \|\partial_{x_3}N_1(\omega(s),\tilde
  \omega(s))\|_{X^p(m)^3})\dd s\\ \nonumber
  \,&\le\, C \int_0^t \frac{e^{-\mu(t-s)}}{a(t{-}s)^{\frac{1}{p}
  -\frac{1}{2}+\frac{|\beta|}{2}}}(\|\omega(s)\|_{\X(m)}\|\nabla 
  \tilde\omega_h(s)\|_{(X(m))^2}\\ \nonumber
  & \hspace{3.8cm} +\|\partial_{x_3}\omega(s)\|_{\X(m)}\|\nabla \tilde
  \omega(s)\|_{\X(m)} +\|\omega(s)\|_{\X(m)}\|\partial_{x_3}\nabla\tilde
  \omega(s)\|_{\X(m)})\dd s\\ \label{est.Phi1}
  \,&\le\, C \int_0^t \frac{e^{-\mu(t-s)} e^{-(\mu+\eta)s}}{a(t{-}s)^{
  \frac{1}{p}-\frac{1}{2}+\frac{|\beta|}{2}} a(s)^\frac{1}{2}}\dd s
  ~\|\omega\|_\V \|\tilde\omega\|_\V \,\le\, C e^{-\mu t}
  \|\omega\|_\V\|\tilde\omega\|_\V~.
\end{align}
In the last inequality, we have used the definition of the norm 
in $\V$ and the fact that $\mu + \eta < 1 + 2\eta$. The bound 
\eqref{est.Phi1} also holds for $t < 1$, and can easily be established
using \eqref{est.prop.linear.1} instead of \eqref{est.prop.linear.3}.

Next, to bound $\partial_{x_3}\Phi_{1,h}(\omega,\tilde\omega)$, we 
recall that $\partial_{x_3}e^{t(L-\alpha\Lambda)} = e^{-t}e^{t(L-\alpha\Lambda)}
\partial_{x_3}$. Applying \eqref{est.prop.linear.1}, we find 
\begin{align*}
  \|\partial_{x_3}\partial_x^\beta &\Phi_{1,h}(\omega,\tilde\omega)(t)\|_{X(m)^2} 
  \,\le\, \int_0^t e^{-(t-s)}\|\partial_x^\beta \{e^{(t-s)(L-\alpha\Lambda)}
  \partial_{x_3}N_1(\omega(s),\tilde\omega(s))\}_h\|_{X(m)^2}\dd s\\
  \,&\le\, C \int_0^t \frac{e^{-(\mu+1)(t-s)}}{a(t{-}s)^{\frac{1}{p}
  -\frac{1}{2}+\frac{|\beta|}{2}}} \|\partial_{x_3} N_1(\omega(s),
  \tilde\omega(s))\|_{X^p(m)^3}\dd s\\
  \,&\le\, C \int_0^t \frac{e^{-(\mu+1)(t-s)}e^{-(\mu+\eta)s}}{a(t{-}s)^{
  \frac{1}{p}-\frac{1}{2}+\frac{|\beta|}{2}}a(s)^\frac{1}{2}}\dd s
  ~\|\omega\|_\V\|\tilde\omega\|_\V  \,\le\, \frac{Ce^{-(\mu+\eta)t}}{a(t)^{
  \frac{1}{p}+\frac{|\beta|}{2}-1}}\,\|\omega\|_\V\|\tilde\omega\|_\V~.
\end{align*}
Similarly, for $k = 0,1$, we can estimate $\partial_{x_3}^k\Phi_{2,h}
(\omega,\tilde\omega)$ as follows:
\begin{align*}
  \|\partial_{x_3}^k\partial_x^\beta &\Phi_{2,h}(\omega,\tilde\omega)
  (t)\|_{X(m)^2} \,\le\, \int_0^t e^{-k(t-s)}\|\partial_x^\beta 
  \{e^{(t-s)(L-\alpha\Lambda)} \partial_{x_3}^k N_2(\omega(s),\tilde
  \omega(s))\}_h\|_{X(m)^2}\dd s\\
  \,&\le\, C \int_0^t \frac{e^{-(\mu+k)(t-s)}}{a(t{-}s)^{\frac{1}{p}
  -\frac{1}{2}+\frac{|\beta|}{2}}} \|\partial_{x_3}^k N_{2}(\omega(s),
  \tilde\omega(s))\|_{X^p(m)^3}\dd s\\
  \,&\le\, C \int_0^t \frac{e^{-(\mu+k)(t-s)}e^{-(\mu+\eta)s}}{a(t{-}s)^{
  \frac{1}{p}-\frac{1}{2}+\frac{|\beta|}{2}}a(s)^\frac{k}{2}}\dd s
  ~\|\omega\|_\V\|\tilde\omega\|_\V \,\le\, \frac{Ce^{-(\mu+k\eta)t}}{
  a(t)^{\frac{1}{p}+\frac{|\beta|}{2}+\frac{k}{2}-\frac{3}{2}}} 
  \|\omega\|_\V\|\tilde\omega\|_\V~.
\end{align*}
Finally, using \eqref{est.prop.linear.2}, we obtain for the vertical 
components of $\Phi_j(\omega,\tilde\omega)$: 
\begin{align*}
  \|\partial_{x_3}^k\partial_x^\beta &\Phi_{j,3}(\omega,
  \tilde\omega)(t)\|_{X(m)} \,\le\, \int_0^t e^{-k(t-s)}\|\partial_x^\beta 
  \{e^{(t-s)(L-\alpha\Lambda)} \partial_{x_3}^k N_j(\omega(s),\tilde
  \omega(s))\}_3\|_{X(m)}\dd s\\
  \,&\le\, C \int_0^t \frac{e^{-(\eta+k)(t-s)}}{a(t{-}s)^{\frac{1}{p}
  -\frac{1}{2}+\frac{|\beta|}{2}}}\|\partial_{x_3}^k N_j(\omega(s),
  \tilde\omega(s))\|_{X^p(m)^3}\dd s\\
  \,&\le\, C \int_0^t \frac{e^{-(\eta+k)(t-s)}e^{-(k+2\eta)s}}{a(t{-}s)^{
  \frac{1}{p}-\frac{1}{2}+\frac{|\beta|}{2}}a(s)^\frac{k}{2}}\dd s
  ~\|\omega\|_\V\|\tilde\omega\|_\V  \,\le\, \frac{Ce^{-(\eta+k)t}}{
  a(t)^{\frac{1}{p}+\frac{|\beta|}{2} +\frac{k}{2}-\frac{3}{2}}} 
  \|\omega\|_\V\|\tilde\omega\|_\V~.
\end{align*}
Summarizing, we have shown that $\Phi$ defined by \eqref{def.Phi}
maps $\V$ into $\V$ and satisfies the following bounds:
\begin{align*}
  \|\Phi(\omega)\|_\V \,&\le\, C_3\sum_{|\beta|\le 1}\|\partial_x^\beta
  \omega_0\|_{\X(m)} + C_4 \|\omega\|_\V^2~, \\
  \|\Phi(\omega)- \Phi(\tilde\omega)\|_\V \,&\le\, C_4(\|\omega\|_\V
  +\|\tilde\omega\|_\V)\|\omega-\tilde\omega\|_\V~, 
\end{align*}
for all $\omega, \tilde \omega \in \V$. If $K = 2C_3\sum_{|\beta|\le 1}
\|\partial_x^\beta\omega_0\|_{\X(m)}$ is sufficiently small, it follows
that $\Phi$ is a strict contraction in the ball $\tilde B_K = \{\omega 
\in \V\,|\, \|\omega\|_\V \le K\}$, hence has a unique fixed point 
there. Denoting by $\omega(t)$ the solution of \eqref{eq.omega.integral}
given by Proposition~\ref{prop.nonlinear.stability}, this implies 
that $t \mapsto \omega(t+T)$ belongs to $\tilde B_K$ if $T > 0$
is sufficiently large. In particular, $\omega(t)$ satisfies 
\eqref{est.thm.main2.1} for some suitable $C > 0$. The proof of 
Theorem~\ref{thm.main2} is now complete. \QED


\section{Appendix}\label{appendix}

\subsection{Proof of Lemma~\ref{lem.reduction}}\label{subsec.lem.reduce}
\label{subsec.reduc}

Let $\chi \in C_0^\infty(\R^2)$ be a cut-off function such that 
$\chi(x_h) = 1$ if $|x_h|\le 1$ and $\chi(x_h) = 0$ if $|x_h| \ge 2$. 
Given $R > 0$, we denote $\chi_R(x_h) = \chi(x_h/R)$, so that 
$|\nabla_h\chi_R(x_h)| \le C/R$. For any $x_3 \in \R$, we define
\[
  f(x_3) \,=\, \int_{\R^2}\tilde\omega_3(x_h,x_3)\dd x_h~, \qquad
  f_R(x_3) \,=\, \int_{\R^2}\tilde\omega_3(x_h,x_3)\chi_R(x_h)\dd x_h~.
\]
Since $\tilde\omega_3\in X(m)$ for some $m > 1$, it is easy to verify 
that $\|f-f_R\|_{L^\infty(\R)}\to 0$ as $R\to \infty$. 
On the other hand, for any test function $\psi \in C_0^\infty(\R)$, 
we have
\begin{equation}\label{eq.hReq}
  \Bigl|\int_{\R}f(x_3)\frac{\D\psi}{\D x_3}(x_3) \dd x_3\Bigr| \,\le\, 
  \Bigl|\int_{\R}f_R(x_3)\frac{\D\psi}{\D x_3}(x_3) \dd x_3\Bigr| 
  + \|f-f_R\|_{L^\infty (\R)} \Bigl\|\frac{\D\psi}{\D x_3}\Bigr\|_{L^1(\R)}~.
\end{equation}
The last term in the right-hand side converges to zero as $R \to 
\infty$. To treat the other term, we observe that
\[
  \int_{\R}f_R(x_3)\frac{\D\psi}{\D x_3}(x_3) \dd x_3 \,=\, 
  \int_{\R^3}\tilde\omega_3(x_h,x_3)\chi_R(x_h)\frac{\D\psi}{\D x_3}
  (x_3) \dd x_3 \,=\, \langle \tilde \omega_3\,,\,\frac{\partial\phi_R}{
  \partial x_3}\rangle~, 
\]
where $\phi_R(x_h,x_3) = \chi_R(x_h)\psi(x_3)$ and $\langle \cdot,\cdot
\rangle$ denotes the duality pairing of $\mathcal{D}'(\R^3)$ and
$C_0^\infty(\R^3)$. Now, since $\nabla\cdot\tilde \omega = 0$ in 
the sense of distributions, we have
\[
  \langle \tilde \omega_3\,,\,\frac{\partial\phi_R}{\partial x_3}
  \rangle \,=\, -\langle \frac{\partial\tilde\omega_3}{\partial x_3}
  \,,\, \phi_R\rangle \,=\, \langle \nabla_h \cdot\tilde \omega_h
  \,,\, \phi_R\rangle \,=\, -\langle \tilde \omega_h
  \,,\, \nabla_h \phi_R\rangle~,
\]
so that
\[
  \int_{\R}f_R(x_3)\frac{\D\psi}{\D x_3}(x_3)\dd x_3 \,=\,  
  -\int_{\R^3} \tilde\omega_h(x_h,x_3)\cdot \nabla_h\chi_R(x_h) 
  \psi(x_3) \dd x_h \dd x_3~.
\]
Using the inclusion $L^2(m) \hookrightarrow L^1(\R^2)$ and the definition 
\eqref{def.X(m)} of the space $X(m)$, we thus find
\[
  \Bigl|\int_{\R}f_R(x_3)\frac{\D\psi}{\D x_3}(x_3)\dd x_3\Bigr|
  \,\le\, \frac{C}{R}\,\|\tilde\omega_h\|_{X(m)^2}\|\psi\|_{L^1(\R)}
  ~\xrightarrow[R \to \infty]{}~ 0~.
\]
Returning to \eqref{eq.hReq}, we conclude that the left-hand 
side vanishes for all $\psi \in C_0^\infty(\R)$, hence $\frac{\D f}{\D 
x_3}=0$ in the sense of distributions. Since $f \in BC(\R)$, it follows
that $f$ is identically constant, which is the desired result. 
\QED

\medskip\noindent{\bf Remark.} If $\omega(x,t)$ is any solution
of \eqref{eq.omega} that is integrable with respect to the 
horizontal variables, we can define
\[
  \phi(x_3,t) \,=\, \int_{\R^2}\omega_3(x_h,x_3,t)\dd x_h~, 
  \qquad x_3 \in \R~, \quad t \ge 0~.
\]
As was observed in \cite{GW2}, this quantity satisfies a remarkably
simple equation
\begin{equation}\label{eq.phi}
  \partial_t \phi(x_3,t) + x_3\partial_{x_3}\phi(x_3,t) \,=\,
  \partial_{x_3}^2\phi(x_3,t)~, 
\end{equation}
which can be solved explicitly. However, if $\omega(\cdot,t) \in
X(m)^3$ for some $m > 1$ with $\nabla\cdot\omega(\cdot,t) = 0$, 
Lemma~\ref{lem.reduction} shows that $\phi(x_3,t)$ does not depend
on $x_3$, and \eqref{eq.phi} then implies that $\phi(x_3,t)$ is 
also independent of $t$. Thus, as was already mentioned, we can 
restrict ourselves to the particular case where $\phi \equiv 0$ 
without loss of generality. Being unaware of this simple observation, 
the authors of \cite{GW2} have stated their stability result in
a seemingly more general form, allowing (apparently) for nontrivial
functions $\phi(x_3,t)$, but thanks to Lemma~\ref{lem.reduction}
(which also holds in the slightly different functional setting of 
\cite{GW2}) the simpler presentation adopted here in Theorem~\ref{thm.main1}
is exactly as general. 

\subsection{Proof of Proposition~\ref{prop.discrete.spectrum}.} 

\label{subsec.discrete}

This final section is devoted to the proof of 
Proposition~\ref{prop.discrete.spectrum}, which shows that 
eigenfunctions of $\cL_{\alpha,h}$ corresponding to eigenvalues 
outside the essential spectrum have a Gaussian decay at infinity.
For the nonlocal operator $\cL_{\alpha,3}$, the same result was 
established in \cite[Lemma 4.5]{GW1} using ODE techniques, but we 
prefer using here a more flexible method based on weighted $L^2$ 
estimates. In fact, we shall consider a more general elliptic
problem of the form
\begin{equation}
  -\LL f + F(x,f,\nabla f) + \lambda f \,=\, h~, \qquad x\in \R^n~,
  \label{eq.prop.appendix.prop.reduction}
\end{equation}
where the unknown is the vector-valued function $f = (f_1,\dots,f_N)^\top$. 
Here and below we denote by $\LL = \Delta+\frac{x}{2}\cdot\nabla 
+\frac{n}{2}$ the analog of operator \eqref{def.LLh} in dimension $n$. 
The data of the problem are the functions $F : \R^n\times \mathbb{C}^N 
\times \mathbb{C}^{nN} \to \mathbb{C}^N$ and $h : \R^n \to \mathbb{C}^N$, 
and the complex number $\lambda$. 

For $m \in [0,\infty]$, we denote by $L^2(m)$, $H^1(m)$ the 
following complex Hilbert spaces on $\R^n$:
\begin{align*}
  L^2(m) \,&=\,\Bigl\{f\in L^2(\R^n,\mathbb{C})~\Big|~\int_{\R^n}|f(x)|^2
  \rho_m(|x|^2)\dd x < \infty\Bigr\}~,\\
  H^1(m) \,&=\, \Bigl\{f\in L^2(m)~\Big|~\partial_{x_j}f\in L^2(m)
  \quad \hbox{for }j = 1,\dots,n\Bigr\}~,
\end{align*}
where $\rho_m$ is the weight function defined by \eqref{def.weight}. 
Our main result is:

\begin{prop}\label{prop.appendix.prop.reduction} Let $m \in [0,\infty)$, 
$\lambda\in \mathbb{C}$, $h\in L^2(\infty)^N$, and assume that $F$ is a 
continuous function satisfying
\begin{equation}
  |F(x,p,Q)| \,\le\, A(x)|p| + B(x)|Q|~, \qquad \hbox{for all}~(x,p,Q)
  \in \R^n\times \mathbb{C}^N\times \mathbb{C}^{nN}~,
  \label{con.prop.appendix.prop.reduction2}
\end{equation}
where $A$ and $B$ are bounded, nonnegative functions such that
\begin{equation}
  \lim_{R\to \infty}\sup_{|x|\ge R} A(x) \,=\, \lim_{R\to
  \infty}\sup_{|x|\ge R} B(x) \,=\, 0~.
  \label{con.prop.appendix.prop.reduction}
\end{equation}
If $\Re \lambda > \frac{n}{4}-\frac{m}{2}$, then any solution $f\in 
H^1(m)^N$ of \eqref{eq.prop.appendix.prop.reduction} satisfies 
$f\in H^1(\infty)^N$.
\end{prop}

\noindent{\bf Proof.} The proof is a simple modification of \cite[Proposition
12]{KM}, which in turn is inspired by a recent work of Fukuizumi and 
Ozawa \cite{FO} where decay estimates are obtained for solutions of
the Haraux-Weissler equation. For $k\ge 1$, $\epsilon>0$, and $\theta
\in [0,m]$, we define the weight functions
\begin{equation}
  \xi_{k,\epsilon}(x) \,=\, e^{\frac{(1-\epsilon)k|x|^2}{4k+|x|^2}}~, \qquad
  \zeta_{\theta}(x) \,=\, (1+|x|^2)^\theta~, \qquad x \in \R^n~.
\end{equation}
Multiplying both sides of \eqref{eq.prop.appendix.prop.reduction} by 
$\zeta_\theta\xi_{k,\epsilon} \bar f$ and integrating by parts the 
real part of the resulting expression, we obtain the identity
\begin{align}
  \int_{\R^n} \zeta_\theta \xi_{k,\epsilon} &|\nabla f|^2 \dd x + \Re
  \int_{\R^n} \bar f \cdot (\nabla  ( \zeta_\theta \xi_{k,\epsilon} ),
  \nabla )f \dd x +  \int_{\R^n}|f|^2 \frac{x}{4} \cdot \nabla 
  (\zeta_\theta \xi_{k,\epsilon} ) \dd x \nonumber \\
  \,&=\,-\Re \int_{\R^n}   \zeta_\theta \xi_{k,\epsilon} \bar f \cdot
  F(x,f(x),\nabla f(x)) \dd x+  \Bigl(\frac{n}{4} -\Re\lambda
  \Bigr)\int_{\R^n}\zeta_\theta \xi_{k,\epsilon} |f|^2 \dd x 
    \label{eq.prop.appendix.prop.reduction2}\\
  &\quad~ + \Re \int_{\R^n}\zeta_\theta \xi_{k,\epsilon}\bar f\cdot h 
  \dd x~.\nonumber
\end{align}
Clearly, 
\begin{equation}\label{eq.nablaxi}
  \nabla\xi_{k,\epsilon}(x) \,=\, \frac{8(1-\epsilon)k^2x}{(4k+|x|^2)^2}
  \,\xi_{k,\epsilon}(x)~, \qquad 
  \nabla\zeta_\theta(x) \,=\, \frac{2\theta x}{1+|x|^2}\,\zeta_\theta(x)~.
\end{equation}
Thus, the second term in the left-hand side of
\eqref{eq.prop.appendix.prop.reduction2} can be written in the following
way:
\begin{align*}
  \Re \int_{\R^n} &\bar f \cdot (\xi_{k,\epsilon} \nabla\zeta_\theta,
  \nabla)f \dd x + \Re \int_{\R^n} \bar f \cdot (\zeta_\theta\nabla 
  \xi_{k,\epsilon},\nabla)f \dd x\\
  \,&=\, -\int_{\R^n} |f|^2\,\nabla\cdot\Bigl(\frac{\theta x\zeta_\theta 
  \xi_{k,\epsilon}}{1+|x|^2}\Bigr)\dd x + \Re \int_{\R^n} \bar f \cdot 
  (\zeta_\theta\nabla \xi_{k,\epsilon},\nabla)f \dd x\\
  \,&=\, -\int_{\R^n} |f|^2\,\xi_{k,\epsilon}\,x \cdot\nabla \frac{\theta 
  \zeta_\theta }{1+|x|^2} \dd x  - \int_{\R^n}|f|^2\frac{\theta 
  \zeta_\theta }{1+|x|^2}\,x \cdot\nabla \xi_{k,\epsilon} \dd x \\
  &\quad~ -n\theta\int_{\R^n} \frac{\zeta_\theta \xi_{k,\epsilon}}{1+|x|^2}
  |f|^2 \dd x + \Re \int_{\R^n} \frac{8(1-\epsilon) k^2\zeta_\theta
  \xi_{k,\epsilon}}{(4k+|x|^2)^2}\,\bar f \cdot (x, \nabla) f \dd x~.
\end{align*}
To bound this quantity from below, we observe that
\[
  \int_{\R^n} |f|^2\,\xi_{k,\epsilon} x \cdot \nabla \frac{\theta 
  \zeta_\theta }{1+|x|^2} \dd x\le 2\theta^2\int_{\R^n} \frac{\zeta_\theta
  \xi_{k,\epsilon}}{1+|x|^2}|f|^2 \dd x~.
\]
Moreover, for each $\eta_1>0$,
\begin{align*}
  - \Re \int_{\R^n} &\frac{8(1-\epsilon) k^2\zeta_\theta \xi_{k,\epsilon}} 
  {(4k+|x|^2)^2}\,\bar f \cdot (x, \nabla) f \dd x \,\le\,
  \int_{\R^n} \frac{2(1-\epsilon)k \zeta_\theta \xi_{k,\epsilon}}{4k+|x|^2}
  |xf| |\nabla f| \dd x \\
  \,&\le\, (1-\eta_1) \int_{\R^n} \zeta_\theta \xi_{k,\epsilon} |\nabla f|^2 
  \dd x + \frac{(1-\epsilon)^2}{1-\eta_1} \int_{\R^n} \frac{ k^2 \zeta_\theta 
  \xi_{k,\epsilon}|xf|^2  }{(4k+|x|^2)^2} \dd x~.
\end{align*}
Thus, using the expression \eqref{eq.nablaxi} of $\nabla \xi_{k,\epsilon}$, 
we find
\begin{align}
  \Re \int_{\R^n} \bar f \cdot (\nabla (\zeta_\theta \xi_{k,\epsilon}),
  \nabla)f &\dd x 
  \,\ge\, -C \int_{\R^n} \frac{\zeta_\theta \xi_{k,\epsilon}}{1+|x|^2}
  |f|^2 \dd x - \int_{\R^n} \frac{8(1-\epsilon)\theta k^2 \zeta_\theta
  \xi_{k,\epsilon}  |xf|^2}{(4k+|x|^2)^2(1+|x|^2)}\dd x 
  \label{eq.prop.appendix.prop.reduction3}\\
  \,&-\, (1-\eta_1) \int_{\R^n} \zeta_\theta \xi_{k,\epsilon} |\nabla f |^2
  dx - \frac{(1-\epsilon)^2}{1-\eta_1} \int_{\R^n} \frac{ k^2 \zeta_\theta
  \xi_{k,\epsilon}|xf|^2  }{(4k+|x|^2)^2} \dd x~, \nonumber
\end{align}
where $C = n\theta+\theta^2$ does not depend on $k$ and $\epsilon$.
We next consider the third term in the left-hand side of
\eqref{eq.prop.appendix.prop.reduction2}, which satisfies
\begin{equation}
  \int_{\R^n}|f|^2 \frac{x}{4}\cdot\nabla (\zeta_\theta\xi_{k,\epsilon}) 
  \dd x \,=\, \frac{\theta}{2} \int_{\R^n} \frac{\zeta_\theta 
  \xi_{k,\epsilon}}{1+|x|^2}\,|xf|^2\dd x + 2(1-\epsilon)\int_{\R^n}
  \frac{k^2 \zeta_\theta \xi_{k,\epsilon} |xf|^2}{(4k+|x|^2)^2} \dd x~.
  \label{eq.prop.appendix.prop.reduction4}
\end{equation}
To estimate the right-hand side of \eqref{eq.prop.appendix.prop.reduction2}, 
we use \eqref{con.prop.appendix.prop.reduction2} and obtain, 
for each $\eta_2>0$,
\begin{align}
  - \Re \int_{\R^n}\zeta_\theta \xi_{k,\epsilon} \bar f \cdot F(x,f(x),
  &\nabla f(x))\dd x \,\le\, \int_{\R^n}\zeta_{\theta}\xi_{k,\epsilon} 
  A|f|^2 \dd x + \int_{\R^n}\zeta_{\theta}\xi_{k,\epsilon}B|f||\nabla f| 
  \dd x \nonumber\\
  \,&\le\, \int_{\R^n}  \zeta_{\theta} \xi_{k,\epsilon}\Bigl(A+\frac{B^2}
  {4\eta_2}\Bigr)|f|^2 \dd x + \eta_2 \int_{\R^n} \zeta_{\theta}\xi_{k,\epsilon}
  |\nabla f|^2 \dd x~. \label{eq.prop.appendix.prop.reduction5} 
\end{align}
Finally, for each $\eta_3>0$, we have
\begin{equation}
  \Re \int_{\R^n}\zeta_\theta \xi_{k,\epsilon} \bar f\cdot h \dd x \,\le\, 
  \eta_3 \int_{\R^n}\zeta_\theta \xi_{k,\epsilon}|f|^2 \dd x + 
  \frac{1}{4\eta_3} \int_{\R^n}\zeta_\theta \xi_{k,\epsilon}|h|^2 \dd x~. 
\label{eq.prop.appendix.prop.reduction6}
\end{equation}
Substituting 
\eqref{eq.prop.appendix.prop.reduction3}--%
\eqref{eq.prop.appendix.prop.reduction6} into 
\eqref{eq.prop.appendix.prop.reduction2}, we arrive at our basic
inequality:
\begin{align}
  (\eta_1-&\eta_2)\int_{\R^n} \zeta_{\theta} \xi_{k,\epsilon} |\nabla f|^2
   \dd x +  \int_{\R^n} \frac{(1-\epsilon) k^2
  \zeta_{\theta} \xi_{k,\epsilon} |x f|^2}{(4k+|x|^2)^2} \Bigl(
  \frac{1-2\eta_1+\epsilon}{1-\eta_1} - \frac{8 \theta}{1+|x|^2}\Bigr)
   \dd x \nonumber \\
  \,&\le\, \int_{\R^n}   \zeta_{\theta} \xi_{k,\epsilon}
  \Bigl\{\Bigl(\frac{C}{1+|x|^2} + \frac{n}{4} - \Re\lambda + A + \frac{B^2}
  {4\eta_2} +\eta_3 -\frac{\theta}{2}\Bigr) |f|^2 + \frac{1}{4\eta_3} 
  |h|^2 \Bigr\} \dd x~. \label{ineq.appendix.prop.reduction}
\end{align} 

To exploit \eqref{ineq.appendix.prop.reduction}, we first take 
$\eta_1=\eta_2=\frac{1}{2}$ and $\theta=m$. Using 
\eqref{con.prop.appendix.prop.reduction} and the assumption that
$\Re\lambda > \frac{n}{4}-\frac{m}{2}$, we see that there exists 
$R > 0$ independent of $k\ge 1$ such that, if $\eta_3>0$ is sufficiently 
small, the following inequality holds: 
\[
  \epsilon(1-\epsilon)\int_{\R^n} \frac{ k^2 \zeta_{\theta}
  \xi_{k,\epsilon} |xf|^2}{(4k+|x|^2)^2} \dd x \le C \int_{|x|\le R} 
  \zeta_{\theta}\xi_{k,\epsilon} |f|^2 \dd x +\frac{1}{4\eta_3}\int_{\R^n}
  \zeta_{\theta}\xi_{k,\epsilon}  |h|^2 \dd x~,
\]
where the constant $C > 0$ is independent of $k\ge 1$. Thus, taking the limit 
$k \to \infty$ and using Fatou's lemma, we obtain 
\[
 \frac{\epsilon(1{-}\epsilon)}{16}\int_{\R^n} (1+|x|^2)^m e^{\frac{1-
 \epsilon}{4}|x|^2} |xf|^2 \dd x \,\le\, C(R) \int_{|x|\le R} |f|^2 \dd x 
 + \frac{1}{4\eta_3}\int_{\R^n} (1+|x|^2)^m e^{\frac{1-\epsilon}{4}|x|^2}
 |h|^2 \dd x~,
\]
which shows that $e^{\frac{1-\epsilon}{8}|x|^2}f \in L^2(\R^2)$ for any 
$\epsilon>0$. Next we choose $\eta_1=\frac{1}{4}$, $\eta_2=\frac{1}{8}$,
$\eta_3=1$, and $\theta=0$ in \eqref{ineq.appendix.prop.reduction}. 
Taking again the limit $k \to \infty$ and using Lebesgue's dominated
convergence theorem, we find
\[
 \frac{1}{8} \int_{\R^n} e^{\frac{1-\epsilon}{4}|x|^2}|\nabla f|^2 \dd x +
  \frac{1{-}\epsilon}{24} \int_{\R^n}  e^{\frac{1-\epsilon}{4}|x|^2} |xf|^2 
  \dd x \le C \int_{\R^n}   e^{\frac{1-\epsilon}{4}|x|^2} |f|^2 \dd x
  +\frac{1}{4}\int_{\R^n} e^{\frac{1-\epsilon}{4}|x|^2}  |h|^2 \dd x~,
\]
where the constant $C > 0$ does not depend on $\epsilon > 0$. This 
inequality shows that 
\[
 \frac{1}{8}\int_{\R^n} \!e^{\frac{1-\epsilon}{4}|x|^2}|\nabla f|^2 \dd x +
 \frac{1{-}\epsilon}{48} \int_{\R^n} \!e^{\frac{1-\epsilon}{4}|x|^2} |xf|^2 
 \dd x \,\le\, C \int_{|x|\le R'} \!e^{\frac{1-\epsilon}{4}|x|^2}|f|^2 
 \dd x + \frac{1}{4}\int_{\R^n} \!e^{\frac{1-\epsilon}{4}|x|^2}|h|^2 \dd x\,,
\]
for some $R'>0$ independent of $\epsilon > 0$. Taking now the limit
$\epsilon\to 0$, we conclude that $f \in H^1(\infty)$, which is the
desired result. \QED

\vspace{0.5cm}

\noindent{\bf Proof of Proposition~\ref{prop.discrete.spectrum}.} 
We consider the eigenvalue equation \eqref{eq.eigenvalue}, which 
can be written in the form
\begin{eqnarray}
  -\LL_h\omega_h + \alpha \Lambda_1\omega_h - \alpha\tilde\Lambda_2
  \omega_h +\Bigl(\lambda +\frac{3}{2}\Bigr)\omega_h \,=\, 0~,
  \label{eq.prop.discrete.spectrum.1}
\end{eqnarray}
where $\LL_h$ is given by \eqref{def.LLh} and the operators 
$\Lambda_1$, $\tilde \Lambda_2$ are defined at the beginning of
Section~\ref{sec.vectorial}. We recall that $|\Lambda_1 \omega_h| \le
|U^G_h| |\nabla_h \omega_h|$ and $|\tilde \Lambda_2 \omega_h| \le
|\nabla_h U^G_h| |\omega_h|$, where the velocity profile $U^G_h$ 
satisfies \eqref{eq.UGbounds}. Assume that $\Re \lambda > -\frac{m}{2}-1$ 
and let $\omega_h \in H^1(m)^2$ be a solution to 
\eqref{eq.prop.discrete.spectrum.1}. Applying 
Proposition~\ref{prop.appendix.prop.reduction} with
$n=N=2$, $F(x,f,\nabla f) = \alpha \Lambda_1 f - \alpha\tilde\Lambda_2 f$, 
and $h=0$, we obtain $\omega_h \in H^1(\infty)^2$. This completes the 
proof of Proposition~\ref{prop.discrete.spectrum}. \QED


\end{document}